\documentclass[12pt]{article}

\pdfoutput=1

\usepackage[T1]{fontenc}

\usepackage[utf8]{inputenc}

\usepackage[english]{babel}

\usepackage[nomath]{lmodern}

\usepackage{amsmath}
\usepackage{amssymb}
\usepackage{amsthm}

\usepackage{stmaryrd}

\usepackage{geometry}

\usepackage{enumitem}

\usepackage{tikz}
	\usetikzlibrary{arrows.meta} 
	\usetikzlibrary{graphs} 
	\usetikzlibrary{calc} 

\usepackage[subrefformat=parens]{subcaption}

\usepackage[pdfusetitle]{hyperref}

\begin{document}
	
	\makeatletter
	

\newcommand*{\cdef}{\newcommand*}

\cdef \credef {%
	\renewcommand*%
}

\credef \cdef {%
	\newcommand*%
}

\cdef \sdef [1]{%
	\cdef #1{}%
	\def #1%
}

\cdef \slet [1]{%
	\cdef #1{}%
	\let #1%
}

\cdef \undef [1]{%
	\let #1 \undefined
}


\slet \Ex \expandafter

\cdef \ExEx {%
	\Ex\Ex\Ex
}

\cdef \ExExEx {%
	\Ex\Ex\Ex\Ex\Ex\Ex\Ex
}

\cdef \ExArg [2]{%
	\Ex#1\Ex{#2}%
}

\cdef \ExExArg [2]{%
	\ExEx#1\ExEx{#2}%
}

\cdef \ExpandedOnce [2]{%
	\begingroup
	\def \x ##1{%
		\endgroup
		#1%
	}%
	\expandafter\x\expandafter{#2}%
}

\cdef \Expanded [1]{%
	\begingroup
	\edef \x {%
		\endgroup
		#1%
	}%
	\x
}

\cdef \Apply [1]{%
	%
	\Expanded{\noexpand #1}%
}


\cdef \IfEmpty [3]{%
	\ifx \undefined#1\undefined
		#2%
	\else
		#3%
	\fi
}


\cdef \IgnoreNext [1]{%
}

\cdef \MacroName [1]{%
	\ExExArg\scantokens{\Ex\IgnoreNext \string#1\noexpand}%
}


\cdef \Show [1]{%
	\expandafter\show\csname #1\endcsname
}

	\undef \int 


\setlist{itemsep = 0pt}
\setlist[enumerate]{leftmargin=*} 
\setlist[enumerate, 1]{label=\upshape (\roman*), ref=(\roman*)}
\setlist[enumerate, 2]{label=\upshape (\alph*), ref=(\alph*)}

\newlist{enumerate_alter}{enumerate}{1}
\setlist[enumerate_alter, 1]{label=\upshape (\alph*), ref=(\alph*)}

\cdef \DescriptionFormat [1]{%
	\normalfont\emph{#1}%
}

\setlist[description]{format=\DescriptionFormat}

\newlist{description_inline}{description*}{1} 
\SetLabelAlign{none}{#1} 
\setlist[description_inline]{format=\DescriptionFormat, align=none}

\cdef \th@slanted {\slshape}	

\theoremstyle{slanted}
\newtheorem{theorem}{Theorem}[section]
\newtheorem{proposition}[theorem]{Proposition}
\newtheorem{lemma}[theorem]{Lemma}
\newtheorem{observation}[theorem]{Observation}
\newtheorem{corollary}[theorem]{Corollary}

\theoremstyle{definition}
\newtheorem{definition}[theorem]{Definition}
\newtheorem{notation}[theorem]{Notation}
\newtheorem{example}[theorem]{Example}
\newtheorem{remark}[theorem]{Remark}
\newtheorem{note}[theorem]{Note}
\newtheorem{problem}[theorem]{Problem}
\newtheorem{question}[theorem]{Question}
\newtheorem*{acknowledgements}{Acknowledgements}


\cdef \@LabelPrefix {}	

\cdef \labelblock [1]{%
	\def \@LabelPrefix {#1}%
	\label{#1}%
}

\cdef \loclabel [1]{%
	\if \@LabelPrefix \relax \relax
		\message{My Warning: No prefix for local label found.}%
	\else
		\label{\@LabelPrefix.#1}%
	\fi
}

\cdef \loclabelblock [1]{%
	\if \@LabelPrefix \relax \relax
		\message{My Warning: No prefix for local label found.}%
	\else
		\ExpandedOnce{\def \@LabelPrefix {##1.#1}}{\@LabelPrefix}%
		\label{\@LabelPrefix}%
	\fi
}

\cdef \locref [1]{%
	\ref*{\@LabelPrefix.#1}%
}

\cdef \itemref [2]{%
	\ref{#1} \ref*{#1.#2}%
}

\cdef \itemonlyref [2]{%
	\ref*{#1.#2}%
}


\cdef \locimpl [2]{%
	{\upshape \locref{#1} $\Longrightarrow$ \locref{#2}}%
}

\cdef \locequiv [2]{%
	{\upshape \locref{#1} $\Longleftrightarrow$ \locref{#2}}%
}

\cdef \ForwardImplication {%
	\textquotedblleft $\Longrightarrow$\textquotedblright%
}

\cdef \BackwardImplication {%
	\textquotedblleft $\Longleftarrow$\textquotedblright%
}

\cdef \A {\mathcal{A}}
\cdef \B {\mathcal{B}}
\cdef \C {\mathcal{C}}
\cdef \F {\mathcal{F}}
\cdef \G {\mathcal{G}}
\credef \P {\mathcal{P}}
\credef \S {\mathcal{S}}
\cdef \T {\mathcal{T}}
\cdef \U {\mathcal{U}}

\cdef \NN {\mathbb{N}}
\cdef \ZZ {\mathbb{Z}}
\cdef \RR {\mathbb{R}}
\cdef \CC {\mathbb{C}}


\cdef \holds {%
	\colon
}

\cdef \st {%
	\colon
}

\cdef \letiff {%
	\; :\Longleftrightarrow \;%
}


\cdef \given {%
	\colon
}

\cdef \set [1]{%
	\{#1\}%
}

\cdef \tuple [1]{%
	\langle #1\rangle
}


\cdef \card [1]{%
	\lvert #1\rvert
}

\cdef \powset [1]{%
	\mathcal{P}(#1)%
}

\cdef \subsets [2]{%
	[#2]^{#1}%
}

\cdef \disunion {%
	\sqcup
}

\cdef \DisUnion {%
	\bigsqcup
}

\cdef \notowns {%
	\mathrel{\m@th\mathpalette\c@ncel\owns}%
}

\cdef \continuum {%
	\mathfrak{c}%
}


\cdef \maps {%
	\colon
}

\cdef \id {%
	\operatorname{id}%
}

\cdef \dom {%
	\operatorname{dom}%
}

\cdef \rng {%
	\operatorname{rng}%
}

\cdef \cod {%
	\operatorname{cod}%
}

\cdef \im [1]{%
	[#1]%
}

\cdef \inv {%
	^{-1}%
}

\cdef \preim [1]{%
	\inv\im{#1}%
}

\cdef \fiber [1]{%
	\inv(#1)%
}

\cdef \restr {%
	\mathop{\upharpoonright}%
}

\cdef \diag {%
	\mathbin\vartriangle
}

\cdef \Diag {%
	\bigtriangleup
}

\cdef \codiag {%
	\mathbin\triangledown
}

\cdef \CoDiag {%
	\bigtriangledown
}


\cdef \homeo {%
	\cong
}

\cdef \clo [2][]{%
	\overline{#2}\IfEmpty{#1}{}{^{#1}}%
}

\cdef \cl {%
	\operatorname{cl}%
}

\cdef \int {%
	\operatorname{int}%
}

\cdef \ro {%
	\operatorname{ro}%
}

\cdef \topsum {%
	\oplus
}

\cdef \TopSum {%
	\sum
}

\cdef \QuotientPart [1]{%
		#1^\textup{q}%
	}

\cdef \InjectivePart [1]{%
	#1^\textup{i}%
}


\cdef \half {%
	\frac{1}{2}%
}

\cdef \mon [1]{%
	#1^\circ
}

\cdef \SG [1]{%
	G^<_{#1}%
}


\cdef \Sep {%
	\mathcal{S}%
}

\cdef \SepF {%
	\Sep_1%
}

\cdef \SepH {%
	\Sep_2%
}

\cdef \SepU {%
	\Sep_{2\half}%
}

\cdef \SepT {%
	\Sep_\textup{\kern0.05em f\kern0.05em}%
}

\cdef \SepZ {%
	\Sep_\textup{c}%
}

\cdef \SymPoints {%
	^\textup{P}%
}

\cdef \Points {%
	^\textup{PP}%
}

\cdef \PointClosed {%
	^\textup{PC}%
}

\tikzset{
	between/.style args={#1 and #2}{
		at = ($(#1)!0.5!(#2)$)
	},
	edge/.style = {
		draw,
		semithick, 
	},
	open/.style = {
		edge,
		circle, 
		inner sep = 0.09cm, 
	}, 
	closed/.style = {
		open,
		fill,
	}, 
	zigzag_correction/.style = {
		x = 1.12cm,
	},
}

	\makeatother

	\title{Tree sums of maximal connected spaces}
	\author{Adam Bartoš \texorpdfstring{
		\footnote{Charles University, Faculty of Mathematics and Physics, Department of Mathematical Analysis.}
		\footnote{E-mail address: \texttt{drekin@gmail.com}.}
	}{}}
	\date{December 2, 2018}
	\maketitle
	
	\vspace{-2em}
	\begin{center} \em
		Dedicated to the memory of Petr Simon, \\ who introduced me to the beautiful realm of general topology.
	\end{center}
	\vspace{0em}
	
	\begin{abstract}
		A topology $\tau$ on a set $X$ is called maximal connected if it is connected, but no strictly finer topology $\tau^* > \tau$ is connected. We consider a construction of so-called tree sums of topological spaces, and we show how this construction preserves maximal connectedness and also related properties of strong connectedness and essential connectedness. 
		
		We also recall the characterization of finitely generated maximal connected spaces and reformulate it in the language of specialization preorder and graphs, from which it is clear that finitely generated maximal connected spaces are precisely $T_\half$-compatible tree sums of copies of the Sierpiński space.
		
		\begin{description}
			\item[Classification:] 
				54A10, 
				54D05, 
				54B17, 
				54D10, 
				54G15. 
			
			\item[Keywords:] maximal connected, strongly connected, essentially connected, tree sum, I-subset, submaximal, nodec, specialization preorder.
		\end{description}
	\end{abstract}
	
	\linespread{1.2}\selectfont

\section{Introduction}
	
	For every fixed set $X$ we may consider the collection of all topologies on $X$. These form a complete lattice $\T(X)$ when ordered by inclusion. For every topological property $\P$ we may consider the subcollection of $\T(X)$ consisting of all topologies having the property $\P$. Then we may consider maximal and minimal elements of this collection.
	A topology $\tau \in \T(X)$ is called \emph{maximal $\P$} if it satisfies $\P$ but no strictly finer topology in $\T(X)$ satisfies $\P$. 
	The property of being \emph{minimal $\P$} is defined dually. Often, the maximality is considered when $\P$ is stable under coarser topologies, and minimality is considered when $\P$ is stable under finer topologies.
	
	Probably the most classical result in this context is the fact that Hausdorff compact spaces are both minimal Hausdorff and maximal compact (but not every minimal Hausdorff space is compact and not every maximal compact space is Hausdorff). References, many other properties, and a general treatment can be found in a paper by Cameron \cite{Cameron_71}.
	There are also \emph{maximal} spaces where the implicit property $\P$ is “having no isolated points”. Let us mention van Douwen's example of countable regular maximal space that can be found in \cite{van_Douwen_93}.
	
	We are interested in the situation where $\P$ means connectedness, i.e.\ in \emph{maximal connected} spaces. These were first considered by Thomas in \cite{Thomas_68}, where he proved among other results that every open connected subspace of a maximal connected space is maximal connected, and also characterized finitely generated maximal connected spaces. There are also related notions of \emph{strongly connected} and \emph{essentially connected} spaces. Following Cameron, for a topological property $\P$ when we consider being maximal $\P$, we say that a topological space is \emph{strongly $\P$} if it admits finer maximal $\P$ topology. Essentially connected are those connected spaces whose every connected expansion has the same connected subsets – these spaces were considered by Guthrie and Stone in \cite{GS_73}.
	
	In this paper we first recall the facts about maximal, strongly, and essentially connected spaces that we use later. Clearly the construction of topological sum does not preserve the properties since it does not preserve connectedness. In the second section we consider another sum-like construction – a \emph{tree sum}. It is a certain quotient of a topological sum – such quotient that it preserves the original spaces as subspaces, glues them only at individual points, and the overall structure of gluing corresponds to a tree graph. First, we systematically treat the properties of tree sums of topological spaces, so we may next show how this construction preserves maximal, strong, and essential connectedness (Theorem~\ref{thm:maximal_connected_tree_sum} and \ref{thm:maximal_connected_tree_sum_equivalence}).
	
	In the third section we revise Thomas' characterization of finitely generated maximal connected spaces. We describe them in the language of specialization preorder and graphs. With this description their structure is crystal clear, they can be easily visualized, and it is evident that they are exactly $T_\half$-compatible tree sums of copies of the Sierpiński space (Corollary~\ref{thm:finitely_generated_maximal_connected_tree_sums}).

	\begin{definition}
		Let $X$ be a set or a topological space. We say that $X$ is \emph{degenerate} if $\card{X} \leq 1$. Otherwise, we say that $X$ is \emph{nondegenerate}.
		
		By a \emph{decomposition} of $X$ we mean an indexed family $\tuple{A_i: i \in I}$ of subsets of $X$ such that $\bigcup_{i \in I} A_i = X$ and $A_i \cap A_j = \emptyset$ for every $i \neq j \in I$. If additionally every set $A_i$ is nonempty, we say that the decomposition is \emph{proper}. Hence, a topological space is connected if and only if it admits no clopen proper decomposition $\tuple{U, V}$.
	\end{definition}
	
	\begin{notation}
		Let $X$ be a set. The order of the lattice $\T(X)$ of all topologies on $X$ is denoted simply by $\leq$. So $\tau \leq \tau^*$ means that $\tau$ is coarser and $\tau^*$ is finer. Also, $\tau < \tau^*$ means that $\tau^*$ is \emph{strictly} finer than $\tau$.
		Additionally, when $\tau \leq \tau^*$ (or $\tau < \tau^*$), we say that the topology $\tau^*$ is an \emph{expansion} of $\tau$ (or a \emph{strict expansion} of $\tau$). In that case we also say that the space $\tuple{X, \tau^*}$ is an expansion of $\tuple{X, \tau}$.
		
		The join operation on $\T(X)$ is denoted by $\vee$ and is extended to all subsystems of $\powset{X}$, so $\A \vee \B$ denotes the topology generated by $\A \cup \B$ for any $\A, \B \subseteq \powset{X}$. Hence, for $\tau$ a topology on $X$ and $\A \subseteq \powset{X}$ the expansion of $\tau$ by $\A$ is denoted by $\tau \vee \A$. For $A \subseteq X$ the expansion $\tau \vee \set{A}$ is called a \emph{simple expansion} of $\tau$.
		
		For a topology $\tau \in \T(X)$ and a set $Y \subseteq X$ the induced subspace topology on $Y$ is denoted by $\tau\restr{Y}$.
	\end{notation}
	
	\begin{definition}
		Recall that a topological space $\tuple{X, \tau}$ or its topology $\tau$ is called
		\begin{itemize}
			\item \emph{maximal connected} if it is connected and has no connected strict expansion;
			\item \emph{strongly connected} if it has a maximal connected expansion;
			\item \emph{essentially connected} if it is connected and every connected expansion has the same connected subsets.
		\end{itemize}
	\end{definition}
	
	\begin{observation} \label{thm:testing_expansions}
		When testing maximal or essential connectedness, it is enough to consider only expansions by finite families. Let $\tuple{X, \tau}$ be a connected topological space and let $\tau^*$ be a connected expansion of $\tau$.
		\begin{enumerate}
			\item If $A \in \tau^* \setminus \tau$, then $\tau' := \tau \vee \set{A}$ is also a connected expansion of $\tau$. Hence, $\tau$ is maximal connected if and only if it has no connected strict simple expansion. Also note that every $\tau'$-open set is of the form $U \cup (A \cap V)$ and also of the form $(U' \cup A) \cap V'$ where $U$, $V$, $U'$, $V'$ are $\tau$-open.
			\item If $C \subseteq X$ is not $\tau^*$-connected, then there are $\tau^*$-open sets $U, V \subseteq X$ such that $\tuple{C \cap U, C \cap V}$ is a proper decomposition of $C$. Hence, for $\tau' := \tau \vee \set{U, V}$ we have that $C$ is $\tau'$-disconnected while $\tau'$ is connected. Therefore, it is enough to test essential connectedness on expansions by two of sets.
		\end{enumerate}
	\end{observation}
	
	The following lemmata provide conditions to test whether a set open in an expansion is open in the original topology as well, and whether a subspace of an expansion is a subspace of the original space as well.
	
	\begin{lemma} \label{thm:originally_open}
		Let $\tuple{X, \tau}$ be a topological space, let $\tau^* = \tau \vee \A$ be an expansion of $\tau$ for some $\A \subseteq \powset{X}$. If a set $U$ is $\tau^*$-open, then the following conditions are equivalent.
		\begin{enumerate}
			\item $U$ is $\tau$-open. 
			\item For every $A \in \A$ there is a $\tau$-open set $V_A$ such that $U \cap A \subseteq V_A \subseteq U \cup A$. 
		\end{enumerate}
		
		\begin{proof}
			For \ForwardImplication{} it is enough to put $V_A := U$. For \BackwardImplication{} note that the set $U$ is of the form $\bigcup_{i \in I} (W_i \cap \bigcap_{j \in J_i} A_{i, j})$ where the sets $W_i$ are $\tau$-open, the sets $A_{i, j}$ are members of $\A$, and the index sets $J_i$ are finite. 
			Consider the function $f$ that maps every indexed family $\tuple{X_{i, j}: i \in I, j \in J_i}$ of subsets of $X$ to the set $\bigcup_{i \in I} (W_i \cap \bigcap_{j \in J_i} X_{i, j})$. Clearly, $f$ is monotone in the sense that for every two families $\tuple{X_{i, j}}$, $\tuple{Y_{i, j}}$ such that $X_{i, j} \subseteq Y_{i, j}$ for every $i, j$ we have $f(\tuple{X_{i, j}}) \subseteq f(\tuple{Y_{i, j}})$. Note that $f(\tuple{U \cap A_{i, j}}) = f(\tuple{U \cup A_{i, j}}) = U$. Therefore, the $\tau$-open set $f(\tuple{V_{A_{i, j}}})$ is also equal to $U$.
		\end{proof}
	\end{lemma}
	
	\begin{corollary} \label{thm:originally_open_corollary}
		Let $\tuple{X, \tau}$ be a topological space, let $\tau^* = \tau \vee \A$ be an expansion of $\tau$ for some $\A \subseteq \powset{X}$. A $\tau^*$-open set $U$ is $\tau$-open if for every $A \in \A$ any of the following conditions holds.
		\begin{enumerate}
			\item $U \cap A$ is $\tau$-open, in particular $U \cap A = \emptyset$. 
			\item $U \cup A$ is $\tau$-open, in particular $U \cup A = X$. 
			\item There is a $\tau$-open set $V$ such that $U \supseteq V \supseteq A$. 
			\item There is a $\tau$-open set $V$ such that $U \subseteq V \subseteq A$. 
		\end{enumerate}
	\end{corollary}
	
	\begin{lemma} \label{thm:unaffected_subspace}
		Let $\tuple{X, \tau}$ be a topological space, let $\tau^* = \tau \vee \A$ be an expansion of $\tau$ for some $\A \subseteq \powset{X}$. If $Y \subseteq X$, then $\tau^*\restr{Y} = \tau\restr{Y}$ if and only if $Y \cap A$ is $\tau$-open in $Y$ for every $A \in \A$, in particular if for every $A \in \A$ we have $Y \cap A = \emptyset$ or $Y \subseteq A$.
		
		\begin{proof}
			\ForwardImplication{} is obvious since $Y \cap A$ is $\tau^*$-open in $Y$ for every $A$. For \BackwardImplication{} let $U \subseteq Y$ be $(\tau^*\restr{Y})$-open, and for every $A \in \A$ let $V_A \subseteq X$ be a $\tau$-open set such that $V_A \cap Y = A \cap Y$. There are $\tau$-open sets $W_i$, finite sets $J_i$, and sets $A_{i, j} \in \A$ for $i \in I$, $j \in I_j$ such that
			$U = \bigcup_{i \in I} (Y \cap W_i \cap \bigcap_{j \in J_i} A_{i, j}) = \bigcup_{i \in I} (Y \cap W_i \cap \bigcap_{j \in J_i} V_{A_{i, j}})$, which is $(\tau\restr{Y})$-open.
		\end{proof}
	\end{lemma}

	\begin{definition} \label{def:weaker_properties}
		Recall that a topological space $X$ is called 
		\begin{itemize}
			\item \emph{submaximal} if every dense subset is open, equivalently if every co-dense subset is closed (and so discrete), equivalently if $\clo{A} \setminus A$ is closed for every $A \subseteq X$;
			\item \emph{nodec} if every nowhere dense subset is closed, equivalently if every nowhere dense subset is discrete, equivalently if $\clo{U} \setminus U$ is discrete for every open $U \subseteq X$;
			\item \emph{$T_\half$} if the singleton $\set{x}$ is open or closed for every $x \in X$.
		\end{itemize}
		Note that co-dense sets are exactly sets of the form $\clo{A} \setminus A$ for $A \subseteq X$, and that closed nowhere dense sets are exactly sets of the form $\clo{U} \setminus U$ for open $U \subseteq X$.
	\end{definition}
	
	Submaximal spaces without isolated points were introduced by Hewitt, who called them MI-spaces. The name submaximal is due to Bourbaki. Many references and an excellent overview can be found in \cite{AC_95}. Nodec spaces were considered by van Douwen in \cite[1.14]{van_Douwen_93}. $T_\half$ spaces were introduced by McSherry in \cite{McSherry_74} under name $T_{ES}$. The name $T_\half$ comes from Levine, who earlier introduced a different but equivalent condition in \cite{Levine_70}.
	
	\begin{proposition} \label{thm:implications}
		All implications in Figure~\ref{fig:implications} hold.
		
		\input{figures/implications.fig}
		
		\begin{proof}
			Every maximal connected space is submaximal since making any dense subset open (or even a filter of dense subsets \cite[Lemma 1]{Anderson_65}) preserves connectedness. Every submaximal space is $T_\half$ since every singleton is either open or co-dense and hence closed. The other implications are clear from definitions.
		\end{proof}
	\end{proposition}

\subsection{Preservation under subspaces}
	
	\begin{observation}	\label{thm:easy_subspace_preservation}
		The properties of being submaximal, nodec, or $T_\half$ are hereditary to all subspaces.
		
		\begin{proof}
			This was proved before, see for example \cite[Proposition 2.1]{AC_95}. It is enough to observe that every co-dense, nowhere dense, or one-point subset of a subspace is co-dense, nowhere dense, or one-point in the original space as well, respectively.
		\end{proof}
	\end{observation}
	
	Clearly, the properties of maximal connectedness, strong connectedness, and essential connectedness can be preserved only by connected subspaces. Thomas proved in \cite[Theorem 3]{Thomas_68} that maximal connectedness is hereditary with respect to open connected subspaces. Later, Guthrie, Reynolds, and Stone proved the same first for closed connected subspaces in \cite[Lemma 2]{GRS_73}, and then using submaximality they observed that every connected subspace of a maximal connected space is open in its closure, and so is itself maximal connected \cite[Theorem 7]{GRS_73}. In \cite[Theorem 1]{GS_73} Guthrie and Stone proved that essential connectedness is hereditary with respect to connected subspaces as well. The core argument of the proofs can be stated as follows.
	
	\begin{lemma} \label{thm:connected_subspace_expansion}
		Let $\tuple{Y, \sigma}$ be a subspace of a connected space $\tuple{X, \tau}$. For every connected expansion $\sigma^* \geq \sigma$ there exists a connected expansion $\tau^* \geq \tau$ such that $\tau^*\restr{Y} = \sigma^*$.
		
		\begin{proof}
			We put $\tau^* := \tau \vee \A$ where $\A := \set{S \cup (X \setminus \clo{Y}): S \in \sigma^*}$. Clearly, $\tau^*$ is an expansion of $\tau$ such that $\tau^*\restr{Y} = \sigma^*$. We need to show that it is connected. $\clo{Y}$ is $\tau^*$-connected since $Y$ is $\tau^*$-connected and $\clo{Y} = \cl_\tau(Y) = \cl_{\tau^*}(Y)$. Let $\tuple{U, V}$ be a $\tau^*$-clopen decomposition of $X$. Without loss of generality $\clo{Y} \subseteq U$. $U$ is $\tau$-open by Corollary~\ref{thm:originally_open_corollary} since $U \cup A = X$ for every $A \in \A$. Let $W$ be the $\tau$-open set $X \setminus \clo{Y}$. $V$ is $\tau$-open by Corollary~\ref{thm:originally_open_corollary} since $V \subseteq W \subseteq A$ for every $A \in \A$. Hence, $\tuple{U, V}$ is a $\tau$-clopen decomposition of $X$, so $U = \emptyset$ or $V = \emptyset$ since $\tau$ is connected.
		\end{proof}
	\end{lemma}
	
	Now, the above-mentioned results on preservation under connected subspaces can be re-proved easily.
	
	\begin{proposition} \labelblock{thm:subspace_preservation}\hfill
		\begin{enumerate}
			\item Every connected subspace of a maximal connected space is maximal connected. \loclabel{maximal_connected}
			\item Every connected subspace of an essentially connected space is essentially connected.\loclabel{essentially_connected}
			\item Every connected subspace of a strongly connected and essentially connected space is both strongly connected and essentially connected. \loclabel{strongly_connected}
		\end{enumerate}
		
		\begin{proof}
			Let $\tuple{X, \tau}$ be a topological space and let $\tuple{Y, \sigma}$ be its connected subspace.
			\begin{enumerate}
				\item Let $\sigma^*$ be a connected expansion of $\sigma$. By Lemma~\ref{thm:connected_subspace_expansion} there is a connected expansion $\tau^* \geq \tau$ such that $\tau^*\restr{Y} = \sigma^*$. Since $\tau$ is maximal connected, we have that $\tau^* = \tau$ and so $\sigma^* = \sigma$.
				\item Let $C \subseteq Y$ be connected and let $\sigma^*$ be a connected expansion of $\sigma$. By Lemma~\ref{thm:connected_subspace_expansion} there is $\tau^*$ a connected expansion of $\tau$ such that $\tau^*\restr{Y} = \sigma^*$. $C$ is $\tau^*$-connected since $\tau$ is essentially connected, and hence $C$ is $\sigma^*$-connected.
				\item Let $\tau^*$ be a maximal connected expansion of $\tau$. Since $\tau$ is essentially connected, $Y$ is $\tau^*$-connected, and hence $\tau^*\restr{Y}$ is a maximal connected expansion of $\sigma$ by \locref{maximal_connected}.
				\qedhere
			\end{enumerate}
		\end{proof}
	\end{proposition}
	
	\begin{example}
		Not every connected subspace of a strongly connected space is strongly connected. By \cite[Theorem 15]{GS_73} no Hausdorff connected space with a dispersion point is strongly connected. Cantor's leaky tent is such a space. Yet, it is a subspace of $\RR^2$, which is strongly connected by \cite[Corollary 5A]{GSW_78} and also by Corollary~\ref{thm:Euclidean_strongly_connected}.
	\end{example}
	
	\begin{observation} \label{thm:strongly_essentially_connected_reals}
		The interval $[0, 1]$ is both strongly connected and essentially connected. Hence, the same holds for the real line $\RR$ and the interval $[0, 1)$.
		
		\begin{proof}
			The fact that $[0, 1]$ is essentially connected was first proved in \cite[Theorem 4.2]{Hildebrand_67} and also follows from \cite[Theorem 10]{GS_73}.
			A maximal connected expansion of $[0, 1]$ was constructed independently in \cite{Simon_78} and \cite{GSW_78}.
			The equivalence of $\RR$, $[0, 1]$, and $[0, 1)$ with respect to having the properties follows from Proposition~\ref{thm:subspace_preservation} and from the fact that $\RR \hookrightarrow [0, 1) \hookrightarrow [0, 1] \hookrightarrow \RR$.
		\end{proof}
	\end{observation}

\begin{section}{Tree sums of topological spaces}

	\begin{definition} \label{def:tree_sum}
		By a \emph{gluing structure} $\G$ we mean an indexed family of topological spaces $\tuple{X_i: i \in I}$ together with an equivalence $\sim$ on $\TopSum_{i \in I} X_i$. These are exactly the data needed to form a \emph{glued sum} $X_\G := \TopSum_{i \in I} X_i / {\sim}$, which is a quotient of the topological sum. We denote the associated canonical maps $X_i \to X_\G$ by $e_{\G, i}$ and the canonical quotient map $\TopSum_{i \in I} X_i \to X_\G$ by $q_\G$.
		
		We define the \emph{set of gluing points} by $S_\G := \set{x \in X_\G: \card{q_\G\fiber{x}} > 1}$, and we define the \emph{gluing graph} $G_\G$ as the quiver (a directed graph allowing multiple edges between a pair of vertices) such that the set of vertices is $I \disunion S_\G$, and $\tuple{s, i, x}$ is a directed edge from $s$ to $i$ if and only if $e_{\G, i}(x) = s$. Even though the edges are directed in order to stress the bipartite nature of the graph, we consider graph notions like connectedness or paths in the corresponding undirected version of the graph if not stated otherwise.
		
		We say that $\G$ \emph{induces a tree sum} if the corresponding gluing graph $G_\G$ is a tree, i.e.\ for every pair of distinct vertices there is a unique path connecting them. In that case, $X_\G$ is called the \emph{tree sum} of $\G$ and the spaces $X_i$ are called \emph{summands}.
		
		Often, when the gluing structure is implied, we just write “$X$ is a glued/tree sum of $\tuple{X_i: i \in I}$”, or “$X := \TopSum_{i \in I} X_i / {\sim}$ is a glued/tree sum” when we want to name the equivalence. In that case, we write $e_{X, i}$, $q_X$, $S_X$, $G_X$ or even $e_i$, $q$, $S$, $G$ (with a short reminder) instead of $e_{\G, i}$, $q_\G$, $S_\G$, $G_\G$, respectively.
	\end{definition}
	
	\begin{remark}
		Despite the lengthy definition above, the notion of tree sum is quite natural. We just glue topological spaces in a way that the spaces are preserved, two spaces may be glued only at one point, and the global structure of connections forms a tree.
	\end{remark}
	
	\begin{remark}
		Because of the connectedness of the gluing graph, all the summands of a tree sum have to be nonempty unless the whole space is empty.
	\end{remark}
	
	\begin{example}
		A wedge sum, that is a space $\TopSum_{i \in I} X_i / {\sim}$ such that one point is chosen in each space $X_i$ and $\sim$ glues these points together, is an example of a tree sum.
	\end{example}
	
	\begin{example}
		The Arens' space, which is the canonical example of a sequential space that is not Fréchet–Urysohn (see \cite[Example~1.6.19]{Engelking_89}), is a certain tree sum of convergent sequences.
	\end{example}
	
	\begin{observation} \label{thm:injective_sum}
		Let $X := \TopSum_{i \in I} X_i / {\sim}$ be a glued sum. All the maps $e_i\maps X_i \to X$ are injective if and only if there is at most one edge between any two vertices in $G_X$.
		
		\begin{proof}
			For $i \in I$ the map $e_i$ is injective if and only if there are no points $x \neq y \in X_i$ such that $e_i(x) = e_i(y) \in S_X$, that is if and only if there are no points $x \neq y \in X_i$ and $s \in S_X$ such that $\tuple{s, i, x}$ and $\tuple{s, i, y}$ are edges in $G_X$.
		\end{proof}
	\end{observation}
	
	\begin{proposition} \label{thm:tree_sum_retractions}
		Let $X$ be a tree sum of spaces $\tuple{X_i: i \in I}$. There exist unique maps $q_i\maps X \to X_i$ for $i \in I$ such that $q_i \circ e_j = \id_{X_{i = j}}$ if $i = j$ and $q_i \circ e_j$ is constant if $i \neq j$. Hence, all the maps $e_i$ are embeddings, all the maps $q_i$ are quotients (even retractions), and all the spaces $X_i$ are retracts of $X$.
		
		\begin{proof}
			By Observation~\ref{thm:injective_sum} all the maps $e_i$ are injective, and hence we may assume $X_i \subseteq X$ (as sets) for every $i \in I$.
			
			Let $i \in I$. The map $q_i$ has to be the identity on $X_i$, and for $j \neq i$ it has to be constant on $X_j$.
			There exists a unique path between $i$ and $j$ in $G_X$, which goes through vertices $i, s_0, i_0, s_1, \ldots, i_{n - 1}, s_n, j$.
			Since $s_n \in X_j$, the constant value of $q_i$ on $X_j$ has to be $q_i(s_n)$, but since $s_n \in X_{i_{n - 1}} \owns s_{n - 1}$, the constant value of $q_i$ on $X_{i_{n - 1}}$ has to be the same and also equal to $q_i(s_{n - 1})$, and so on. Hence, the constant value on $X_j$ is $q_i(s_0) = s_0$. This is a consistent definition of $q_i$, and it is the only possible.
			
			Since the topology on $X$ is inductively generated by the maps $e_i$, and $q_i \circ e_j$ is continuous for every $i, j \in I$, all the maps $q_i$ are continuous. And since $q_i \circ e_i = \id_{X_i}$, we have that $q_i$ is a quotient map and $e_i$ is an embedding.
		\end{proof}
	\end{proposition}

\begin{subsection}{Internal characterization}
	
	We dedicate a few following paragraphs to an internal characterization of tree sums. Similarly to other sum-like constructions it makes sense to ask if a given topological space is an “inner” tree sum of some of its subspaces.
	
	\begin{definition}
		Let $f\maps X \to Y$ be a continuous map between topological spaces $X$, $Y$. Recall, that there is an induced equivalence relation $\sim_f$ on $X$: $x \sim_f y$ if and only if $f(x) = f(y)$. There is an induced quotient map $\QuotientPart{f}\maps X \to X / {\sim_f}$ and an induced map $\InjectivePart{f}\maps X / {\sim_f} \to Y$ such that $\InjectivePart{f} \circ \QuotientPart{f} = f$. We call $\QuotientPart{f}$ and $\InjectivePart{f}$ the \emph{quotient part} of $f$ and the \emph{injective part} of $f$, respectively.
		
		Let $\tuple{f_i\maps X_i \to Y}_{i \in I}$ be a family of continuous maps. Recall, that there is a canonical map $\CoDiag_{i \in I} f_i\maps \TopSum_{i \in I} X_i \to Y$ called the \emph{codiagonal sum} and defined by the equalities $(\CoDiag_{i \in I} f_i) \circ e_j = f_j$ for $j \in I$ where $e_j\maps X_j \to \TopSum_{i \in I} X_i$ are the canonical embeddings.
	\end{definition}
	
	\begin{definition} \label{def:inner_tree_sum}
		Let $X$ be a topological space and $\F := \tuple{X_i: i \in I}$ a family of its subspaces. The \emph{gluing structure induced by $ \F$} is $\G := \tuple{\F, \sim_f}$ where $f := \CoDiag_{i \in I} e_i\maps \TopSum_{i \in I} X_i \to X$ and $e_i\maps X_i \to X$ are the embeddings for $i \in I$. 
		The gluing structure $\G$ induces the glued sum $X_\G = \TopSum_{i \in I} X_i / {\sim_f}$ by Definition~\ref{def:tree_sum}.
		We say that $\tuple{X, \F}$ is an \emph{inner tree sum}, or that $X$ is an \emph{inner tree sum of $\F$}, if $\G$ induces a tree sum and $\InjectivePart{f}\maps X_\G \to X$ is a homeomorphism.
		
		The family $\F$ also induces a set $S_\F := \set{x \in X: \card{\set{i \in I: x \in X_i}} \geq 2}$ and a graph $G_\F$ on $I \disunion S_\F$ where $\tuple{s, i}$ is an edge from $s \in S_\F$ to $i \in I$ if and only if $s \in X_i$. Note that $S_\F$ and $G_\F$ are canonically isomorphic to $S_\G$ and $G_\G$, respectively. We often identify $\F$ with $\G$, and we write $S_X$ and $G_X$ instead of $S_\F$ and $G_\F$ when the family $\F$ is implied.
	\end{definition}
	
	\begin{remark}
		We need the outer tree sum to construct bigger spaces from summands, but when a bigger space is already constructed, we usually assume that the summands are subspaces of the sum, and we switch to the inner view.
		
		Even though we define the inner tree sum so that the connection with the outer tree sum is clear, the following characterization is easier to work with.
	\end{remark}
	
	\begin{proposition} \label{thm:inner_tree_sum}
		Let $X$ be a topological space and $\F := \tuple{X_i: i \in I}$ a family of its subspaces. $X$ is an inner tree sum of $\F$ if and only if the following conditions hold.
		\begin{enumerate}
			\item $\bigcup_{i \in I} X_i = X$.
			\item $X$ is inductively generated by the family $\F$.
			\item The undirected version of $G_\F$ is a tree.
		\end{enumerate}
		
		\begin{proof}
			Let $e_i\maps X_i \to X$ and $e'_i\maps X_i \to \sum_{j \in I} X_j$ be the canonical embeddings for every $i \in I$. Let us consider the map $f := \CoDiag_{i \in I} e_i\maps \TopSum_{i \in I} X_i \to X$. Clearly, $\InjectivePart{f}$ is bijective iff $f$ is surjective iff $\bigcup_{i \in I} X_i = X$. Note that $\TopSum_{i \in I} X_i / {\sim_f}$ is inductively generated by the family $\tuple{\QuotientPart{f} \circ e'_i: i \in I}$. By the universal property of inductive generation, a bijective $\InjectivePart{f}$ is a homeomorphism if and only if $X$ is inductively generated by the family $\tuple{\InjectivePart{f} \circ \QuotientPart{f} \circ e'_i = f \circ e'_i = e_i: i \in I}$. Finally, $\TopSum_{i \in I} X_i / {\sim_f}$ is a tree sum if and only if $G_\F$ is a tree since $G_\F$ is canonically isomorphic to the gluing graph of $\TopSum_{i \in I} X_i / {\sim_f}$.
		\end{proof}
	\end{proposition}
	
	\begin{observation} \label{thm:irrelevant_singletons}
		Let $X$ be a topological space and $\F := \tuple{X_i: i \in I}$ a family of its subspaces. If $\tuple{x_j: j \in J}$ is a family of points in $\bigcup_{i \in I} X_i$ and $\F' := \F \disunion \tuple{\set{x_j}: j \in J}$, then $X$ is a tree sum of $\F$ if and only if $X$ is a tree sum of $\F'$ (i.e.\ one-point spaces in gluing structures are essentially irrelevant).
		
		\begin{proof}
			Clearly, $\bigcup\rng(\F) = \bigcup\rng(\F')$. Also, $X$ is inductively generated by $\F$ if and only if it is inductively generated by $\F'$ since every member of $\F'$ is contained in a member of $\F$. Finally, $G_\F$ is a tree if and only if $G_{\F'}$ is a tree. We have $S_{\F'} = S_\F \cup \set{x_j: j \in J}$, and the vertices of $G_{\F'}$ are $I \disunion S_{\F'} \disunion J$. Also, $G_\F$ is the subgraph of $G_{\F'}$ induced by $I \disunion S_\F$. For every $s \in S_{\F'} \setminus S_\F$ there is exactly one $i \in I$ such that $s \in X_i$, and the graph $G_{\F'}$ adds $s$ as a new gluing vertex and $\tuple{s, i}$ as a new edge. The graph $G_{\F'}$ also adds every $j \in J$ as a new vertex and $\tuple{x_j, j}$ as a new edge. These changes clearly do not affect, whether the graph is a tree.
		\end{proof}
	\end{observation}
	
\end{subsection}

\begin{subsection}{Tree subsums and branches}
	
	\begin{definition}
		Let $X$ be a tree sum of a family of its subspaces $\tuple{X_i: i \in I}$ and let $Y \subseteq X$. We often use the following notation.
		\begin{itemize}
			\item $I_Y := \set{i \in I: Y \cap X_i \neq \emptyset}$ and $I_x := I_{\set{x}}$ for $x \in X$.
			\item $S_Y := S_X \cap Y$.
			\item $G_Y$ denotes the subgraph of $G_X$ induced by $I_Y \disunion S_Y$.
			\item $\F_Y := \tuple{Y \cap X_i: i \in I_Y}$.
		\end{itemize}
		We say that $Y$ is a \emph{tree subsum} of $X$ if it is an inner tree sum of the family $\F_Y$. Note that $S_Y = S_{\F_Y}$ and $G_Y = G_{\F_Y}$, so the notation is consistent with Definition~\ref{def:inner_tree_sum}.
	\end{definition}
	
	\begin{proposition} \labelblock{thm:tree_subsum}
		Let $X$ be a tree sum of spaces $\tuple{X_i: i \in I}$ and $Y \subseteq X$. The following conditions are equivalent.
		\begin{enumerate}
			\item $Y$ is a tree subsum of $X$.  \loclabel{definition}
			\item $G_Y$ is connected, i.e.\ it is a subtree of $G_X$.  \loclabel{connected_graph}
			\item $q_i\im{Y} = Y \cap X_i$ for every $i \in I_Y$.  \loclabel{projections}
		\end{enumerate}
		
		\begin{proof}
			Let $Y_i$ denote $Y \cap X_i$ for every $i \in I$.
			\begin{description}
				\item[\locimpl{definition}{connected_graph}] is trivial.
				
				\item[\locimpl{connected_graph}{definition}.]
					By Proposition~\ref{thm:inner_tree_sum} it remains to show that $Y$ is inductively generated by the subspaces $\tuple{Y_i: i \in I_Y}$. Let $U \subseteq Y$ be such that for every $i \in I_Y$ the set $U \cap Y_i$ is open in $Y_i$, i.e.\ there is $U_i \subseteq X_i$ open in $X_i$ such that $U_i \cap Y_i = U \cap Y_i$.
					
					Let us put $W := \bigcup_{i \in I_Y} U_i$. We have $W \cap Y = \bigcup_{i \in I_Y} U_i \cap Y = \bigcup_{i \in I_Y} U_i \cap Y_i = \bigcup_{i \in I_Y} U \cap Y_i = U$. But $W$ does not have to be open in $X$. For every $i \in I_U$ we consider the set $V_i := \bigcup\set{q_i\fiber{s}: s \in S_{U_i} \setminus Y_i}$. These are all points in the summands of $X$ attached to $X_i$ via some gluing point in $U_i \setminus Y$. For every $x \in V_i$ there is $j \in I \setminus \set{i}$ such that $x \in X_j$ and a unique path from $j$ to $i$ in $G_X$. This path goes through some $s \in S_{U_i} \setminus Y_i$. Since $G_Y$ is connected and $i \in G_Y$ and $s \notin G_Y$, we have $j \notin G_Y$ and $x \notin Y$. Therefore, $V_i \cap Y = \emptyset$.
					
					Let us put $W' := W \cup \bigcup_{i \in I_U} V_i$. We have $W' \cap Y = W \cap Y = U$. To show that $W'$ is open in $X$ it is enough to observe that $W' \cap X_i$ is open in $X_i$ for every $i \in I$. Let $i \in I_U$; from the definition of $V_i$, we have $V_i \cap X_j \in \set{\emptyset, X_j}$ for every $j \neq i$, and $V_i \cap X_i = S_{U_i} \setminus Y_i \subseteq U_i$. Therefore, $W' \cap X_i \in \set{\emptyset, X_i}$ for every $i \in I \setminus I_U$, and $W' \cap X_i = U_i$ for every $i \in I_U$.
				
				\item[\locimpl{connected_graph}{projections}.]
					Let $i \in I_Y$. We have $q_i\im{Y} = Y_i \cup \set{s \in S_{X_i}: Y \cap q_i\fiber{s} \setminus \set{s} \neq \emptyset}$. Let $s \in S_{X_i}$ and $x \in q_i\fiber{s} \setminus \set{s}$. There is $j \in I$ such that $x \in X_j$ and the path from $i$ to $j$ in $G_X$ goes through $s$. We have $i \in G_Y$, so if $x \in Y$, then $j \in G_Y$ and $s \in G_Y$ by the connectedness of $G_Y$, and hence $s \in Y_i$.
				
				\item[\locimpl{projections}{connected_graph}.]
					Let $i \in G_Y$, $s \in S_{X_i}$, and $j \in G_Y$ such that the path from $i$ to $j$ in $G_X$ goes through $s$. It is enough to show that $s \in Y$. We have that $\set{s} = q_i\im{Y_j} \subseteq q_i\im{Y} = Y_i \subseteq Y$.
				\qedhere
			\end{description}
		\end{proof}
	\end{proposition}
	
	\begin{proposition} \labelblock{thm:subfamily_of_subspaces}
		Let $X$ be a tree sum of spaces $\tuple{X_i: i \in I}$. Let $\F := \tuple{Y_i: i \in I_\F}$ be a family such that $I_\F \subseteq I$ and $Y_i \subseteq X_i$ for every $i \in I_\F$. Let $Y := \bigcup_{i \in I_\F} Y_i$.
		\begin{enumerate}
			\item If $G_\F$ is connected, then $Y \cap X_i = Y_i$ for every $i \in I_\F$.  \loclabel{consistent_intersections}
			\item If $G_\F$ is connected, then $\card{Y \cap X_i} \leq 1$ for every $i \in I \setminus I_\F$.  \loclabel{degenerate_intersections}
			\item If $Y$ is a tree sum of $\F$, then $Y$ is a tree subsum of $X$.  \loclabel{subsum}
		\end{enumerate}
		
		\begin{proof} \hfill
			\begin{enumerate}
				\item Let $i \in I_\F$. Clearly, $Y \cap X_i = Y_i \cup (S_Y \cap X_i)$. For every $s \in S_Y \cap X_i \setminus Y_i$ there is $j \in I_\F \setminus \set{i}$ such that $s \in Y_j$. Since $G_\F$ is connected and it is an induced subgraph of $G_X$, which is a tree, the path $\tuple{i, s, j}$ in $G_X$ is a path in $G_\F$ as well. Hence, $s \in S_\F$ and $s \in Y_i$.
				
				\item Let $i \in I \setminus I_\F$. Suppose that $s \neq s' \in Y \cap X_i$. Let $j, j' \in I_\F$ be such that $s \in Y_j \cap X_i$ and $s' \in Y_{j'} \cap X_i$. We have that $j \neq j'$ since otherwise $\tuple{i, s, j}$ and $\tuple{i, s', j}$ would be two different paths in $G_X$. Hence, we have a path $\tuple{j, s, i, s', j'}$ in $G_X$. But since $G_\F$ is a connected subgraph of $G_X$, there is another path from $j$ to $j'$. That is a contradiction since $G_X$ is a tree.
				
				\item We are comparing the families $\F$ and $\F_Y = \tuple{Y \cap X_i: i \in I_Y}$. We may assume that every $Y_i \neq \emptyset$, otherwise we would have $\F = \tuple{\emptyset}$, $Y = \emptyset$, $\F_Y = \tuple{}$, and the claim would hold. By that assumption, $I_\F \subseteq I_Y$. Since $G_\F$ is a tree, we have $\F = \F_Y\restr{I_\F}$ by \locref{consistent_intersections} and $\card{Y \cap X_i} = 1$ for $i \in I_Y \setminus I_\F$ by \locref{degenerate_intersections}. Therefore, we may use Observation~\ref{thm:irrelevant_singletons}, and $Y$ is a tree sum of $\F_Y$ since it is a tree sum of $\F$.
				\qedhere
			\end{enumerate}
		\end{proof}
	\end{proposition}
	
	\begin{lemma} \label{thm:coarser_family}
		Let $X$ be a tree sum of spaces $\tuple{X_i: i \in I}$. Let $\F := \tuple{Y_j: j \in J}$ be a family of subspaces of $X$. If for every $i \in I$ there is $j \in J$ such that $X_i \subseteq Y_j$, then $\bigcup_{j \in J} Y_j = X$ and $X$ is inductively generated by $\F$. Therefore, $X$ is a tree sum of $\F$ if and only if $G_\F$ is a tree.
		
		\begin{proof}
			We use Proposition~\ref{thm:inner_tree_sum}. Clearly, we have $\bigcup_{j \in J} Y_j \supseteq \bigcup_{i \in I} X_i = X$. Also, if a set $U \subseteq X$ is such that $U \cap Y_j$ is open in $Y_j$ for every $j \in J$, then $U \cap X_i$ is open in $X_i$ for every $i \in I$, and hence $U$ is open in $X$, and hence $X$ is inductively generated by $\F$. The conclusion follows again from Proposition~\ref{thm:inner_tree_sum}.
		\end{proof}
	\end{lemma}
	
	\begin{definition}
		Let $X$ be a tree sum of spaces $\tuple{X_i: i \in I}$. For every $x \in X$ we define its \emph{branches} $\tuple{B_{x, i}: i \in I_x}$ by the formula $B_{x, i} := \set{x} \cup q_i\preim{X_i \setminus \set{x}}$.
		That is, for $x \in S_X$ we have $B_{x, i} = \bigcup\set{X_j: j \in J_{x, i}}$ where $J_{x, i}$ is the set of all indices $j \in I$ such that in the path from $x$ to $j$ in $G_X$ the edge from $x$ goes to $i$. 
		
		Note that if $x \in S_X$ we have $\card{I_x} \geq 2$ and $B_{x, i} \cap B_{x, j} = \set{x}$ for every $i \neq j \in I_x$, whereas if $x \in X \setminus S_X$ there is only one $i$ in $I_x$ and $B_{x, i} = X$.
	\end{definition}
	
	\begin{observation} \label{thm:branches}
		Let $X$ be a tree sum of spaces $\tuple{X_i: i \in I}$, let $x \in X$, and let $\B_x := \tuple{B_{x, i}: i \in I_x}$ be the enumeration of branches at $x$. We have that every $B_{x, i}$ is a tree subsum of $X$ and that $X$ is a tree sum of $\B_x$.
		
		\begin{proof}
			For $x \in X \setminus S_X$ this is clear, so let $x \in S_X$. For every $j \in I$ such that $B_{x, i} \cap X_j \neq \emptyset$ there is a path $\tuple{x, i_0, s_0, \ldots, j_n = j}$ in $G_X$. Either $i_0 = i$ and $X_{j_k} \subseteq B_{x, i}$ for every $k \leq n$, or $i_0 = j$. In both cases the path lies in $G_{B_{x, i}}$, and hence the graph is connected and $B_{x, i}$ is a tree subsum of $X$ by Proposition~\ref{thm:tree_subsum}. For every $i \neq j \in I_x$ we have $B_{x, i} \cap B_{x, j} = \set{x}$, and hence $G_{\B_x}$ is a tree. Since every summand $X_j$ lies in some $B_{x, i}$, we have that $X$ is a tree sum of $\B_x$ by Lemma~\ref{thm:coarser_family}.
		\end{proof}
	\end{observation}
	
\end{subsection}

\begin{subsection}{Separation of tree sums}
	
	Now we will introduce an alternative description of standard separation axioms that is based on existence of continuous maps into special topological spaces. We do this in order to prove preservation of separation axioms for tree sums in a uniform and concise way.
	
	\begin{definition}
		We say that $\Sep$ is a \emph{monotone separation scheme} if $\Sep = \tuple{Y_\Sep, \leq_\Sep}$ where $Y_\Sep$ is a topological space containing the points $0$ and $1$, and $\leq_\Sep$ is a linear order on $Y_\Sep$ such that $0$ is the minimum, $1$ is the maximum, and for every $y \in Y_\Sep$ we have $Y_\Sep / (\leftarrow, y] \homeo [y, \rightarrow)$ via the obvious canonical map (this last condition is motivated by Observation \ref{thm:monotone_separation}).
		
		Let $X$ be a topological space. We say that a pair $A, B \subseteq X$ is $\S$-separated if there is a continuous function $f\maps X \to Y_\Sep$ such that $f\im{A} \subseteq \set{0}$ and $f\im{B} \subseteq \set{1}$. We also say that 
		\begin{itemize}
			\item $X$ is $\Sep\SymPoints$-separated if for every $x \neq y \in X$ either $\set{x}, \set{y}$ or $\set{y}, \set{x}$ is $\Sep$-separated,
			\item $X$ is $\Sep\Points$-separated if every pair of distinct points of $X$ is $\Sep$-separated,
			\item $X$ is $\Sep\PointClosed$-separated if every point and every closed set not containing that point are $\Sep$-separated.
		\end{itemize}
		We consider the following monotone separation schemes:
		\begin{itemize}
			\item $\SepF$ is the Sierpiński space on $\set{0 < 1}$ with isolated point $1$.
			\item $\SepH$ is the space $\set{0 < \half < 1}$ with the topology generated by the two singletons $\set{0}$, $\set{1}$.
			\item $\SepU$ is the space $\set{0 < \bar{0} < \half < \bar{1} < 1}$ with the topology generated by the sets $\set{0}$, $\set{0, \bar{0}, \half}$, $\set{\half, \bar{1}, 1}$, $\set{1}$.
			\item $\SepT$ is $[0, 1]$ where the order is inherited from $\RR$.
			\item $\SepZ$ is the discrete space $\set{0 < 1}$.
		\end{itemize}
	\end{definition}
	
	\begin{observation} \label{thm:monotone_separation}
		Let $\Sep$ be a monotone separation scheme, $X$ a topological space. If $A, B \subseteq X$ are $\Sep$-separated, then for every $y \in Y_\Sep$ there is a continuous map $f\maps X \to Y_\Sep$ such that $f\im{A} \subseteq \set{y}$, $f\im{B} \subseteq \set{1}$.
		
		\begin{proof}
			Let $y \in Y_\Sep$, let $q\maps Y_\Sep \to [y, \rightarrow)$ be the quotient map induced by the canonical homeomorphism $Y_\Sep / (\leftarrow, y] \homeo [y, \rightarrow)$. If $f_0\maps X \to Y_\Sep$ $\Sep$-separates $A, B$, then for $f := q \circ f_0\maps X \to [y, \rightarrow) \subseteq Y_\Sep$ we have $f\im{A} \subseteq \set{q(0)} = \set{y}$, $f\im{B} \subseteq \set{q(1)} = \set{1}$.
		\end{proof}
	\end{observation}
	
	\begin{observation} \label{thm:separation_axioms}
		Let $X$ be a topological space.
		\begin{itemize}
			\item $X$ is $T_0$ if and only if it is $\SepF\SymPoints$-separated.
			\item $X$ is symmetric (i.e.\ for every $x \neq y \in X$, if there is open $U_x$ such that $x \in U_x \notowns y$, then there is open $U_y$ such that $x \notin U_y \owns y$) if and only if for every point $x$ disjoint from a closed set $F$ there is an open set $U$ such that $x \notin U \supseteq F$, that is if and only if $X$ is $\SepF\PointClosed$-separated.
			\item $X$ is $T_1$ if and only if it is $\SepF\Points$-separated.
			\item $X$ is $T_2$ or Hausdorff if and only if it is $\SepH\Points$-separated.
			\item $X$ is $T_{2\half}$ or Urysohn (i.e.\ for every $x \neq y \in X$ there are open sets $U_x \owns x$ and $U_y \owns y$ such that $\clo{U_x} \cap \clo{U_y} = \emptyset$) if and only if it is $\SepU\Points$-separated.
			\item $X$ is functionally $T_2$ (i.e.\ for every $x \neq y \in X$ there is a continuous function $f\maps X \to [0, 1]$ such that $f(x) = 0$ and $f(y) = 1$) if and only if it is $\SepT\Points$-separated.
			\item $X$ is totally separated (i.e.\ for every $x \neq y \in X$ there is a clopen set $U \subseteq X$ such that $x \in U \notowns y$) if and only if it is $\SepZ\Points$-separated.
			\item $X$ is regular if and only if it is $\SepH\PointClosed$-separated.
			\item $X$ is completely regular if and only if it is $\SepT\PointClosed$-separated.
			\item $X$ is zero-dimensional if and only if it is $\SepZ\PointClosed$-separated.
		\end{itemize}
	\end{observation}
	
	\begin{observation} \label{thm:hereditary_separation}
		For every monotone separation scheme $\Sep$ the properties of being $\Sep\SymPoints$-separated, $\Sep\Points$-separated, and $\Sep\PointClosed$-separated are hereditary.
		
		\begin{proof}
			Let $X \subseteq Y$ be topological spaces. If $x$, $y$ are distinct points of $X$, then they are distinct points of $Y$. If $x$ is a point not in a closed set $F$ in $X$, then $x \notin \cl_Y(F)$. Hence, we can move the situation to $Y$. If $f$ $\Sep$-separates the corresponding sets in $Y$, then $f\restr{X}$ $\Sep$-separates the corresponding sets in $X$.
		\end{proof}
	\end{observation}
	
	\begin{proposition}
		Let $X$ be a tree sum of spaces $\tuple{X_i: i \in I}$ and let $\Sep$ be a monotone separation scheme.
		The space $X$ is $\Sep\SymPoints$-, $\Sep\Points$-, or $\Sep\PointClosed$-separated if and only if all the spaces $X_i$ are $\Sep\SymPoints$-, $\Sep\Points$-, or $\Sep\PointClosed$-separated, respectively.
		
		\begin{proof}
			\ForwardImplication{} follows from Observation \ref{thm:hereditary_separation}. \BackwardImplication{} for $\Sep\SymPoints$ and $\Sep\Points$ follows from the fact that for every two points $x \neq y \in X$ there is $i \in I$ such that $q_i(x) \neq q_i(y)$ since if a continuous map $f\maps X_i \to Y_\Sep$ $\Sep$-separates $\set{q_i(x)}$, $\set{q_i(y)}$, then $f \circ q_i\maps X \to Y_\Sep$ $\Sep$-separates $\set{x}$, $\set{y}$.
			
			So let $x \neq y \in X$. If $x \notin S_X$, let $i \in I$ be the index such that $x \in X_i$. If $y \in X_i$, then $q_i(y) = y \neq x = q_i(x)$; if $y \notin X_i$, then $q_i(y) \in S_X \notowns x = q_i(x)$. If $y \notin S_X$, we proceed symmetrically. If $x, y \in S_X$, consider the only path in $G_X$ from $x$ to $y$, going through vertices $x, i_0, s_0, \ldots, i_n, y$. Then we have $q_{i_0}(x) = x \neq s_0 = q_{i_0}(y)$.
			
			To prove \BackwardImplication{} for $\Sep\PointClosed$ let $x \in X$, $x \notin F \subseteq X$ closed. If $x \in S_X$, we put $S' := S_X$ and $G' := G_X$; if $x \notin S_X$, we put $S' := S_X \cup \set{x}$ and define $G'$ as the graph on $I \disunion S'$ extending $G_X$ with the edge $\tuple{x, i_0, x}$ where $i_0 \in I$ is the index such that $x \in X_{i_0}$. We also define a strict partial order $<$ on $G'$: $a < b$ if and only if the path from $x$ to $a$ is a strict initial segment of the path from $x$ to $b$. Basically, we are just rooting the tree at $x$ in order to perform an inductive construction.
			
			We will define continuous maps $f_i\maps X_i \to Y_\Sep$ for $i \in I$ and values $y_s \in Y_\Sep$ for $s \in S'$. We define $y_x := 0$. If $s \in S'$ is a $<$-successor of $i \in I$, we define $y_s := f_i(s)$. If $i \in I$ is a $<$-successor of $s \in S'$, we define $f_i$ as a continuous map such that $f_i(s) = y_s$ and $f_i\im{F \cap X_i} \subseteq \set{1}$, which exists by Observation \ref{thm:monotone_separation}. By the construction, $f := \bigcup_{i \in I} f_i \maps X \to Y_\Sep$ is continuous and $\Sep$-separates $\set{x}$, $F$.
		\end{proof}
	\end{proposition}
	
	\begin{corollary} \label{thm:tree_sum_separation}
		Let $X$ be a tree sum of spaces $\tuple{X_i: i \in I}$. The space $X$ is separated if and only if all the spaces $X_i$ are separated with “separated” meaning $T_0$, symmetric, $T_1$, $T_2$, $T_{2\half}$, functionally $T_2$, totally separated, regular, completely regular, or zero-dimensional.
	\end{corollary}
	
	\begin{definition}
		Let $X$ be a topological space. We say that $A \subseteq X$ is a \emph{$T_\half$-subset} of $X$ if every point of $A$ is closed or isolated in $X$. Equivalently, $X$ is $T_\half$ at every point of $A$.
	\end{definition}
	
	In order to take care of the $T_\half$ separation axiom we need to introduce the following condition.
	
	\begin{definition}
		Let $X$ be a tree sum of spaces $\tuple{X_i: i \in I}$. We say that the corresponding gluing is \emph{$T_\half$-compatible} if we never glue a non-isolated closed point to a non-closed isolated point, i.e.\ there are no $s \in S_X$, $i, j \in I_s$ such that $s$ is non-isolated and closed in $X_i$ and non-closed isolated in $X_j$. Note that if the spaces $X_i$ are $T_\half$, this is equivalent to never gluing a non-closed point to a non-isolated point.
	\end{definition}
	
	\begin{proposition} \labelblock{thm:T_half_tree_sum}
		Let $X$ be a tree sum of spaces $\tuple{X_i: i \in I}$.
		\begin{enumerate}
			\item $S_X$ is a $T_\half$-subset of $X$ if and only if every $S_X \cap X_i$ is a $T_\half$-subset of $X_i$ and the gluing is $T_\half$-compatible.  \loclabel{gluing_set}
			\item The space $X$ is $T_\half$ if and only if all spaces $X_i$ are $T_\half$ and the gluing is $T_\half$-compatible.  \loclabel{whole_space}
		\end{enumerate}
		
		\begin{proof} \hfill
			\begin{enumerate}
				\item Clearly, if $S_X$ is a $T_\half$-subset of $X$, then $S_X \cap X_i$ is a $T_\half$-subset of $X_i$ for every $i \in I$. Also, under this condition the gluing is $T_\half$-compatible if and only if every $s \in S_X$ is closed in every $X_i$ or isolated in every $X_i$ for $i \in I_s$. In other words, if and only if $S_X$ is a $T_\half$-subset of $S_X$.
				
				\item The space $X$ is $T_\half$ if and only if both $S_X$ and $X \setminus S_X$ are $T_\half$-subsets of $X$. The same holds for spaces $X_i$. It is enough to use \locref{gluing_set} and observe that a point in $X \setminus S_X$ is closed or isolated in $X$ if and only if it is so in the space $X_i$ that contains it.
				\qedhere
			\end{enumerate}
		\end{proof}
	\end{proposition}
	
\end{subsection}

\begin{subsection}{Neighborhood-related properties and I-subsets}
	
	Let us start with two lemmata for building neighborhoods in tree sums.
	
	\begin{lemma} \label{thm:piecewise_open}
		Let $X$ be a tree sum of spaces $\tuple{X_i: i \in I}$, let $U_i \subseteq X_i$ for every $i \in I$, and let $U := \bigcup_{i \in I} U_i$. If every $U_i$ is open in the corresponding $X_i$ and every $s \in S_U$ is either isolated or $s \in \bigcap_{i \in I_s} U_i$, then $U$ is open in $X$.
		
		\begin{proof}
			It is enough to show that every $U \cap X_i$ is open in the corresponding $X_i$. Clearly, $U \cap X_i = U_i \cup S_{U \cap X_i}$. By our assumptions, every gluing point in $U \cap X_i$ is either isolated or already contained in $U_i$, and hence $U \cap X_i$ is open in $X_i$.
		\end{proof}
	\end{lemma}
	
	\begin{lemma} \label{thm:star_neighborhood}
		Let $X$ be a tree sum of spaces $\tuple{X_i: i \in I}$ and $x \in X$. Let $S'$ denote the set of all non-isolated gluing points. If $\tuple{U_i: i \in I_x}$ is a family such that for every $i \in I_x$ we have $x \in U_i \subseteq X_i$ and $U_i$ is open in $X_i$, and $V \subseteq X$ is open such that $x \in V$ and $V \cap S' \subseteq \set{x}$, then $W := V \cap \bigcup_{i \in I_x} U_i$ is an open neighborhood of $x$ in $X$.
		
		\begin{proof}
			The claim follows from Lemma~\ref{thm:piecewise_open} applied to the family $\tuple{V \cap U_i: i \in I}$ where we additionally put $U_i := \emptyset$ for $i \in I \setminus I_x$.
		\end{proof}
	\end{lemma}
	
	Let us define the notion of I-subset that naturally occurs in several following propositions.
	
	\begin{definition}
		Let $X$ be a topological space. We say that $A \subseteq X$ is an \emph{\textup{I}-subset} of $X$ if it is a union of an open discrete subset and a closed discrete subset of $X$. Equivalently, the points of $A$ that are not isolated in $X$ form a closed discrete subset of $X$. The name is derived from the related concept of I-space introduced in \cite[Definition 1.4]{AC_95}.
	\end{definition}
	
	\begin{definition}
		Let $X$ be a tree sum of spaces $\tuple{X_i: i \in I}$. We say that the corresponding gluing is \emph{\textup{I}-compatible} if we never glue a non-isolated closed point to an isolated point, i.e.\ there are no $s \in S_X$, $i, j \in I_s$ such that $s$ is non-isolated and closed in $X_i$ and isolated in $X_j$. Note that if the spaces $X_i$ are $T_\half$, this is equivalent to never gluing an isolated point to a non-isolated point.
	\end{definition}
	
	Note that the “I” in the notions of I-subset, I-space, and I-compatibility is a constant symbol referring to “isolated” rather than a mathematical variable $I$. In particular, it is not related to the index set $I$ that is sometimes present.
	
	\begin{observation} \label{thm:T_half/I-subset}
		In a topological space every I-subset is a $T_\half$-subset.
	\end{observation}
	
	\begin{observation} \label{thm:T_half/I-compatibility}
		Let $X$ be a tree sum of spaces $\tuple{X_i: i \in I}$. If the gluing is I-compatible, then it is $T_\half$-compatible. If, additionally, no space $X_i$ contains a clopen point, then the other implication holds as well.
		
		\begin{proof}
			The claim follows from the definition. If there are no clopen points in spaces $X_i$, then every isolated point is non-closed.
		\end{proof}
	\end{observation}
	
	\begin{observation} \labelblock{thm:inductive_discreteness}
		Let $X$ be a topological space inductively generated by a family of its subspaces $\tuple{X_i: i \in I}$. Let $A \subseteq X$ and let $A_i := A \cap X_i$ for every $i \in I$.
		\begin{enumerate}
			\item $A$ is closed discrete if and only if every $A_i$ is closed discrete in the corresponding $X_i$.
				\loclabel{closed}
			\item $A$ is open discrete if and only if every $A_i$ is open discrete in the corresponding $X_i$.
				\loclabel{open}
		\end{enumerate}
		
		\begin{proof}
			Clearly, if $A$ is closed discrete in $X$, then so is every $A_i$ in $X_i$. For the other implication let every $A_i$ be closed discrete in $X_i$. For every $B \subseteq A$ and $i \in I$ we have that $B \cap X_i$ is closed in $X_i$ since $B \cap X_i \subseteq A_i$ and $A_i$ is closed discrete. Therefore, every $B \subseteq A$ is closed in $X$, and hence $A$ is closed discrete in $X$. The proof of \locref{open} is analogous.
		\end{proof}
	\end{observation}
	
	\begin{proposition} \label{thm:gluing_set_discreteness}
		Let $X$ be a tree sum of spaces $\tuple{X_i: i \in I}$.
		\begin{enumerate}
			\item The set $S_X$ is discrete if and only if $S_X \cap X_i$ is discrete for every $i \in I$.
			\item If $S_X$ is an I-subset of $X$, then $S_X \cap X_i$ is an I-subset of $X_i$ for every $i \in I$ and the gluing is $T_\half$-compatible.
			\item If $S_X \cap X_i$ is an I-subset of $X_i$ for every $i \in I$ and the gluing is I-compatible, then $S_X$ is an I-subset of $X$.
		\end{enumerate}
		
		\begin{proof} \hfill
			\begin{enumerate}
				\item Clearly, if $S_X$ is discrete, then $S_X \cap X_i$ is discrete for every $i \in I$. For the other implication let $s \in S_X$. Since $S_X \cap X_i$ is discrete for every $i \in I$, then for every $i \in I_s$ there is $U_i$ open in $X_i$ such that $U_i \cap S_X = \set{s}$. Consider $U := \bigcup_{i \in I_s} U_i$. Then $U \cap S_X = \set{s}$ and $U$ is open since $U \cap X_i = U_i$ if $i \in I_s$, $\emptyset$ otherwise.
				
				\item Clearly, $S_X \cap X_i$ is an I-subset of $X_i$ for every $i \in I$. The rest follows from Observation~\ref{thm:T_half/I-subset} and Proposition~\itemref{thm:T_half_tree_sum}{gluing_set}.
				
				\item Let $S'$ be the set of all points of $S_X$ not isolated in $X$ and let $S'_i$ be the set of all points of $S_X \cap X_i$ not isolated in $X_i$ for every $i \in I$. For every $s \in S'$ there is $i_s \in I_s$ such that $s$ is not isolated in $X_{i_s}$. The point $s$ is also closed in $X_{i_s}$ since $S'_{i_s}$ is closed discrete. Since the gluing is I-compatible, $s$ is not isolated in any $X_i$ for $i \in I$. Therefore, $S' \cap X_i \subseteq S'_i$ for every $i \in I$, and since every $S'_i$ is closed discrete in $X_i$, the set $S'$ is closed discrete in $X$ by Observation~\ref{thm:inductive_discreteness}.
				\qedhere
			\end{enumerate}
		\end{proof}
	\end{proposition}
	
	The following examples show that the claims in Proposition~\ref{thm:gluing_set_discreteness} are sharp.
	
	\begin{example}
		Let $X$ be a wedge sum of spaces $\tuple{X_i: i \in I}$, i.e.\ $S_X = \set{x}$ for some $x \in X$. If $x$ is a closed point in every $X_i$, then $S_X$ is clearly an I-subset of $X$. On the other hand, if additionally $x$ is clopen in some but not all spaces $X_i$, then we glued an isolated point to a closed non-isolated point, so the gluing is not I-compatible.
	\end{example}
	
	\begin{example}
		Let us consider a tree sum $X = \sum_{n \leq \omega} X_n / {\sim}$ where $X_n$ for $n < \omega$ is the Sierpiński space on $\set{0, 1}$ with isolated point $1$, the space $X_\omega$ is the convergent sequence $\omega + 1$, and we glue $\tuple{\omega, n} \sim \tuple{n, 0}$ for every $n < \omega$, i.e.\ we glue the $n$-th member of the sequence with the non-isolated point of the corresponding Sierpiński space.
		
		We have that the gluing is $T_\half$-compatible, and so $X$ is $T_\half$ by Proposition~\ref{thm:T_half_tree_sum}. We also have that $S_X \cap X_n$ is an I-subset of $X_n$ for every $n \leq \omega$, but $S_X$ is not an I-subset of $X$ – it contains no isolated point of $X$ and it is discrete but not closed.
	\end{example}
	
	Now we use the condition of $S_X$ being an I-subset of $X$ as an assumption.
	
	\begin{proposition}
		Let $X$ be a tree sum of spaces $\tuple{X_i: i \in I}$. If $S_X$ is an I-subset of $X$, then $X$ is hereditarily inductively generated by the inclusions of the spaces $X_i$.
		
		\begin{proof}
			Let $U \subseteq A \subseteq X$. Put $A_i := A \cap X_i$ for $i \in I$. We want to show that if $U \cap A_i$ is open in $A_i$ for every $i \in I$, then $U$ is open in $A$. Let $x \in U$. For every $i \in I_x$ there is $U_i$ open in $X_i$ such that $U_i \cap A_i = U \cap A_i$. Let $S'$ be the set of all non-isolated gluing points. Since $S'$ is closed discrete, there is an open set $V \subseteq X$ containing $x$ such that $V \cap S' \subseteq \set{x}$. By Lemma~\ref{thm:star_neighborhood} $W := V \cap \bigcup_{i \in I_x} U_i$ is open in $X$, and we have $x \in W \cap A \subseteq \bigcup_{i \in I_x} U_i \cap A = \bigcup_{i \in I_x} U_i \cap A_i = \bigcup_{i \in I_x} U \cap A_i \subseteq U$.
		\end{proof}
	\end{proposition}
	
	\begin{proposition} \label{thm:tree_sum_expansion}
		Let $\tuple{X, \tau} := \TopSum_{i \in I} \tuple{X_i, \tau_i} / {\sim}$ be a tree sum, let $\A \subseteq \powset{X}$. We put $\tau^* := \tau \vee \A$, $\tau_i^* := \tau_i \vee \set{A \cap X_i: A \in \A}$. If we have that
		\begin{enumerate}
			\item $S_X$ is an I-subset of $\tuple{X, \tau}$,
			\item for every $x \in S_X$ there is a $\tau$-open set $G_x$ such that $\set{A \in \A: x \in A \nsupseteq G_x}$ is finite;
		\end{enumerate}
		then $\tuple{X, \tau^*} = \TopSum_{i \in I} \tuple{X_i, \tau_i^*} / {\sim}$, i.e.\ such expansion of a tree sum is a tree sum of the corresponding expansions.
		
		\begin{proof}
			Clearly, all the maps $e_i\maps \tuple{X_i, \tau_i^*} \to \tuple{X, \tau^*}$ are continuous, and hence we have that $\id_X\maps \TopSum_{i \in I} \tuple{X_i, \tau_i^*} / {\sim} \to \tuple{X, \tau^*}$ is continuous by the inductive generation. To prove the equality it is enough to show that $\tau^*$ is inductively generated by maps $e_i\maps \tuple{X_i, \tau_i^*} \to X$. So let $U \subseteq X$ be such that $U \cap X_i$ is $\tau_i^*$-open for every $i \in I$. We will show that $U$ is $\tau^*$-open.
			
			Let $x \in U \setminus S'$ where $S'$ denotes the set of all non-isolated gluing points. If $x$ is an isolated gluing point, then we are done, otherwise let $i$ be the only $i \in I$ such that $x \in X_i$. Let $U_i$ be a $\tau_i$-open set and $\B \subseteq \A$ a finite family such that $x \in U_i \cap \bigcap\B \subseteq U$. We have that $U_i \setminus S'$ is $\tau$-open since it is $\tau_i$-open and for every $j \neq i$ it holds that $(U_i \setminus S') \cap X_j$ is either empty or an isolated gluing point connecting $X_i$ with $X_j$. Therefore, $W_x := U_i \cap \bigcap\B \setminus S'$ is a $\tau^*$-neighborhood of $x$ in $U$.
			
			Let $x \in U \cap S'$. We put $B := \bigcap\set{A \in \A: x \in A \nsupseteq G_x}$, which is $\tau^*$-open since the set is finite. For every $i \in I_x$ there is an $\tau_i$-open set $U_i \subseteq G_x$ such that $x \in U_i \cap B \subseteq U$. There is also a $\tau$-open set $V$ such that $V \cap S' = \set{x}$. By Lemma~\ref{thm:star_neighborhood} $\bigcup_{i \in I_x} U_i \cap V$ is a $\tau$-neighborhood of $x$, so $W_x := \bigcup_{i \in I_x} U_i \cap V \cap B$ is a $\tau^*$-neighborhood of $x$ in $U$.
		\end{proof}
	\end{proposition}
\end{subsection}

\begin{subsection}{Connectedness-related properties}
	
	Now we focus on connectedness-related properties of tree sums.
	
	\begin{proposition} \label{thm:connected_tree_sum}
		A tree sum $X$ of spaces $\tuple{X_i: i \in I}$ is connected if and only if all the spaces $X_i$ are connected.
		
		\begin{proof}
			\ForwardImplication. By Proposition~\ref{thm:tree_sum_retractions} all the spaces $X_i$ are quotients of the connected space $X$.
			\BackwardImplication. Every component of connectedness contains all spaces $X_i$ that it intersects. Hence, it contains whole $X$ because every two spaces $X_i$, $X_j$ are connected via a path in $G_X$.
		\end{proof}
	\end{proposition}
	
	\begin{observation} \label{thm:closed_or_isolated_point}
		Let $X$ be a tree sum of spaces $\tuple{X_i: i \in I}$. Let $x \in X$ and let $\tuple{B_i: i \in I_x}$ be the branches at $x$. If $x$ is closed or isolated in $X$, then $\tuple{B_i \setminus \set{x}: i \in I_x}$ is a clopen decomposition of $X \setminus \set{x}$.
		
		\begin{proof}
			If $x$ is closed, then every $B_i \setminus \set{x}$ is open in $X$ since the only space $X_j$ such that $(B_i \setminus \set{x}) \cap X_j \notin \set{\emptyset, X_j}$ is $X_i$ where the intersection is $X_i \setminus \set{x}$, which is open. Similarly, if $x$ is isolated, then every $B_i \setminus \set{x}$ is closed in $X$, and hence clopen in $X \setminus \set{x}$.
		\end{proof}
	\end{observation}
	
	\begin{observation} \labelblock{thm:disconnected_co-point}
		Let $X$ be a topological space, let $x \in X$, and let $\tuple{X_i: i \in I}$ be a clopen decomposition of $X \setminus \set{x}$. One of the following situations happens.
		\begin{enumerate}
			\item The point $x$ is closed in $X$ and every $X_i$ is open in $X$.  \loclabel{closed}
			\item The point $x$ is isolated in $X$ and every $X_i$ is closed in $X$.  \loclabel{isolated}
			\item There is $i \in I$ such that $x$ is neither closed nor isolated in $X_i$ while $x$ is clopen in $X_j$ for every $j \in I \setminus \set{i}$.  \loclabel{neither}
		\end{enumerate}
		
		\begin{proof}
			Since $\tuple{X_i: i \in I}$ is a clopen decomposition of $X \setminus \set{x}$, there are sets $\tuple{U_i: i \in I}$ open in $X$ such that for every $i \in I$ we have $U_i \setminus \set{x} = X_i$. If we may choose $U_i = X_i$ for every $i \in I$, we are in situation \locref{closed}. Otherwise, there is $i \in I$ such that $U_i = X_i \cup \set{x}$ and $X_i$ is not open in $X$. If there is $j \in I \setminus \set{i}$ such that we may choose $U_j = X_j \cup \set{x}$, we are in situation \locref{isolated} since $\set{x} = U_i \cap U_j$ and $X \setminus X_k = \bigcup\set{U_l: l \in I \setminus \set{k}}$ for every $k \in I$. If there is no such $j$, then $U_j = X_j$ for every $j \in I \setminus \set{i}$, the point $x$ is not isolated in $X$, and $\tuple{U_i,\, U_j: j \in I \setminus \set{i}}$ is a clopen decomposition of $X$. Hence, we are in situation \locref{neither}.
		\end{proof}
	\end{observation}
	
	\begin{proposition} \label{thm:connected_in_tree_sum}
		Let $X$ be a tree sum of spaces $\tuple{X_i: i \in I}$ such that every gluing point is closed or isolated, i.e. $S_X$ is a $T_\half$-subset of $X$. Let $C \subseteq X$ and $C_i := C \cap X_i$ for every $i \in I_C$. The set $C$ is connected if and only if every $C_i$ is connected and $G_C$ is connected. That is, connected subspaces of $X$ are exactly tree subsums of connected subspaces.
		
		\begin{proof}
			Suppose that $C$ is connected. Let $s \in S_X \setminus C$ and let $\tuple{B_i: i \in I_s}$ be the branches of $X$ at $s$. Since $s$ is closed or isolated and because of Observation~\ref{thm:branches} and \ref{thm:closed_or_isolated_point}, it follows from Observation~\ref{thm:disconnected_co-point} that every $B_i \setminus \set{s}$ is clopen in $X \setminus \set{s}$, and hence $C \subseteq B_i \setminus \set{s}$ for some $i \in I_s$. Therefore, $G_C$ is a connected graph.
		
			If $G_C$ is connected, then $C$ is a tree sum of $\tuple{C_i: i \in I_C}$ by Proposition~\ref{thm:tree_subsum}, and the claim follows from Proposition~\ref{thm:connected_tree_sum} applied to $C$ and $\tuple{C_i: i \in I_C}$.
		\end{proof}
	\end{proposition}
	
\end{subsection}

\begin{subsection}{Maximal connectedness of tree sums}
	
	Now we finally use the machinery built in the previous sections to prove the theorems about maximal connectedness in tree sums of topological spaces.
	
	\begin{theorem} \labelblock{thm:maximal_connected_tree_sum}
		Let $X$ be a tree sum of spaces $\tuple{X_i: i \in I}$ such that the set of all non-isolated gluing points is closed discrete, i.e.\ $S_X$ is an I-subset of $X$.
		\begin{enumerate}
			\item If the spaces $X_i$ are maximal connected, then $X$ is maximal connected.  \loclabel{maximal_connected}
			\item If the spaces $X_i$ are strongly connected, then $X$ is strongly connected.  \loclabel{strongly_connected}
			\item If the spaces $X_i$ are essentially connected, then $X$ is essentially connected.  \loclabel{essentially_connected}
		\end{enumerate}
		
		\begin{proof} Let $\tau$ be the topology on $X$, $\tau_i$ the topology on $X_i$ for every $i \in I$.
			\begin{enumerate}
				\item By Proposition~\ref{thm:connected_tree_sum} $X$ is connected. Let $A \subseteq X$ be non-$\tau$-open. Consider $\tau^* := \tau \vee \set{A}$ and $\tau_i^* := \tau_i \vee \set{A \cap X_i}$ for $i \in I$. Since $A$ is not $\tau$-open, there is $i \in I$ such that $A \cap X_i$ is not $\tau_i$-open, and hence $\tau_i^*$ is disconnected since $\tau_i$ is maximal connected. By Proposition~\ref{thm:tree_sum_expansion} $\tuple{X, \tau^*}$ is a tree sum of the spaces $\tuple{X_i, \tau_i^*}$. Therefore, it is disconnected by Proposition~\ref{thm:connected_tree_sum}.
				
				\item Let $\tau_i^*$ be a maximal connected expansion of $\tau_i$ for every $i \in I$, let $\tuple{X, \tau^*}$ be the corresponding tree sum. Clearly, $\tau^*$ is an expansion of $\tau$. Since $S_X$ is an I-subset of $\tuple{X, \tau}$, it is an I-subset of $\tuple{X, \tau^*}$ as well, and hence by \locref{maximal_connected} $\tau^*$ is maximal connected.
				
				\item Again, $X$ is connected by Proposition~\ref{thm:connected_tree_sum}. By Observation~\ref{thm:testing_expansions} it is enough to test essential connectedness only on expansions by finite families. Let $C$ be a connected subset of $\tuple{X, \tau}$, let $\tau^*$ be a connected expansion of $\tau$ by a finite family, and let $\tau_i^* := \tau^*\restr{X_i}$ for $i \in I$. $\tuple{X, \tau^*}$ is a tree sum of spaces $\tuple{X_i, \tau_i^*}$ by Proposition~\ref{thm:tree_sum_expansion}. By Proposition~\ref{thm:connected_in_tree_sum} every $\tau_i^*$ is connected and every $C_i := C \cap X_i$ is $\tau_i$-connected. By the essential connectedness every $C_i$ is $\tau_i^*$-connected and hence $\tau^*$-connected. Therefore, $C$ is $\tau^*$-connected again by Proposition~\ref{thm:connected_in_tree_sum}.
				\qedhere
			\end{enumerate}
		\end{proof}
	\end{theorem}
	
	\begin{corollary} \label{thm:Euclidean_strongly_connected}
		The following spaces are strongly connected: $\RR^\kappa$ (or generally any real topological vector space), $[0, 1]^\kappa$, $[0, 1)^\kappa$ for $\kappa \geq 1$, spheres, and many others for which the technique described in the proof can be adapted.
		
		\begin{proof}
			Let us consider the space $\RR^\kappa$. Let $\set{L_i: i \in I}$ be the family of all lines through the origin in $\RR^\kappa$, and let $\tau$ be the topology on $\RR^\kappa$ inductively generated by the lines (this idea is due to \cite[Corollary 5A]{GSW_78}). Clearly, $\tau$ refines the standard topology on $\RR^\kappa$. By Proposition \ref{thm:inner_tree_sum} $\tuple{\RR^\kappa, \tau}$ is a tree sum of the lines $L_i$. By Theorem~\ref{thm:strongly_essentially_connected_reals} the lines are strongly connected, and hence $\tuple{\RR^\kappa, \tau}$ is strongly connected by Theorem~\ref{thm:maximal_connected_tree_sum}.
			
			The situation with the other spaces is analogous. The general idea is to cut a space so one gets a tree sum of real intervals that refines the original topology. For example a sphere is cut into meridians glued together at one pole with the other pole attached to one of the meridians.
		\end{proof}
	\end{corollary}
	
	\begin{question}
		Is every CW complex strongly connected?
	\end{question}
	
	\begin{observation}
		Let $\tuple{X, \tau}$ be a topological space, $\tau^*$ a connected expansion of $\tau$, $x \in X$. We denote the family of all $\tau$-connected components of $X \setminus \set{x}$ by $\C_x$ and the family of all $\tau^*$-connected components of $X \setminus \set{x}$ by $\C^*_x$. If $\tuple{X, \tau}$ is essentially connected, then $\C_x = \C^*_x$. Hence, if $\tuple{X, \tau}$ is a topological realization of a tree graph, then any maximal connected expansion $\tau^*$ still possesses the structure of the graph: the vertices of degree $\neq 2$ can be easily identified and the edges remain connected by essential connectedness.
	\end{observation}
	
	Let us focus on the question, whether the assumption of $S_X$ being an I-subset of $X$ is necessary in Theorem~\ref{thm:maximal_connected_tree_sum}.
	
	\begin{example}
		Not every tree sum of copies of the Sierpiński space is maximal connected. Let $X_1$ be the Sierpiński space on $\set{0, 1}$ with isolated point $1$ and $X_2$ the Sierpiński space on $\set{1, 2}$ with isolated point $2$. Consider $X := (X_1 \topsum X_2) / {\sim}$ where $\sim$ glues the points $1$ together. The specialization order (Definition~\ref{def:specialization_preorder}) on $X$ is $0 < 1 < 2$ and the gluing is not $T_\half$-compatible, and hence $X$ is not maximal connected since it is even not $T_\half$.
		
		Also, the set $C := \set{0, 2} \subseteq X$ is connected, but the graph $G_C$ is not connected. This shows that the assumption about closed or isolated gluing points in Proposition~\ref{thm:connected_in_tree_sum} is necessary.
	\end{example}
	
	\begin{example}
		Not every tree sum of maximal connected spaces is maximal connected or even essentially connected. Let $\tuple{[0, 1]_x: x \in [0, 1]}$ be copies of the real interval $[0, 1]$ with a maximal connected expansion of the standard topology. Consider a comb-like space $\tuple{X, \tau} := \TopSum_{x \in [0, 1]} [0, 1]_x / {\sim}$ where $\sim$ glues together points $\tuple{0, x} \sim \tuple{x, 1}$ for $x > 0$. $\tuple{X, \tau}$ is a tree sum of the maximal connected intervals, but it is not maximal connected itself.
		
		Let $A := \set{\tuple{0, 0}} \cup \bigcup_{x > 0} [0, 1)_x$, and $\tau^* := \tau \vee \set{A}$. Clearly, $\tau^*$ is a strict expansion of $\tau$ since $A \cap [0, 1]_0 = \set{\tuple{0, 0}}$, which is not $\tau$-open in $[0, 1]_0$. By Lemma~\ref{thm:unaffected_subspace}, the set $X \setminus \set{\tuple{0, 0}}$ is $\tau^*$-connected since it is $\tau$-connected and $(X \setminus \set{\tuple{0, 0}}) \cap A = \bigcup_{x > 0} [0, 1)_x$, which is $\tau$-open. We also have that $\tuple{0, 0}$ is in $\tau^*$-closure of $\bigcup_{x > 0} [0, 1)_x$. Together, $\tau^*$ is still connected. But since $[0, 1]_0$ becomes disconnected in $\tau^*$, $\tau$ is not even essentially connected.
	\end{example}
	
	\begin{question}
		Is the space $X$ form the previous example strongly connected? Is every tree sum of maximal connected spaces strongly connected?
	\end{question}
	
	The fact that the space from the previous example is not even essentially connected is not a coincidence as the following observation shows.
	
	\begin{observation} \label{thm:essential_implies_maximal}
		Let $X$ be a topological space whose topology is inductively generated by a family of maximal connected subspaces $\tuple{X_i: i \in I}$. If $X$ is essentially connected, then it is maximal connected.
		
		\begin{proof}
			Let $\tau$ be the topology on $X$. If $X$ is not maximal connected, then there is a non-open set $A \subseteq X$ such that $\tau^* := \tau \vee \set{A}$ is still connected. By the inductive generation there is $i \in I$ such that $A \cap X_i$ is not open in $X_i$, and hence $X_i$ is $\tau$-connected but not $\tau^*$-connected, so $\tuple{X, \tau}$ is not essentially connected.
		\end{proof}
	\end{observation}
	
	\begin{proposition} \label{thm:closed_discrete_gluing_points}
		Let $X$ be a tree sum of spaces $\tuple{X_i: i \in I}$ and let $S$ be the set of all points $s \in S_X$ such that there are $i \neq j \in I$ such that $\set{s}$ is nowhere dense in both $X_i$ and $X_j$, i.e.\ $S$ is the set of all points that are nowhere dense in at least two summands. If $X$ is nodec, then $S$ is closed discrete.
		
		\begin{proof}
			By Observation~\ref{thm:inductive_discreteness} it is enough to show that $S \cap X_i$ is closed discrete in $X_i$ for every $i \in I$. Let $i \in I$ and for every $s \in S \cap X_i$ let $\tuple{B_{s, j}: j \in I_s}$ be the enumeration of branches of $X$ at $s$. By the definition of $S$ there is $j \in I_s \setminus \set{i}$ such that $\set{s}$ is nowhere dense in $X_j$. Let $U_s := B_{s, j} \setminus \set{s}$.
			
			\emph{$U_s$ is open in $X$.} Since $\set{s}$ is nowhere dense in $X_j$ and hence in $B_{s, j}$, it is closed in $B_{s, j}$, which is nodec as a subspace of $X$ (Proposition~\ref{thm:easy_subspace_preservation}). Hence, $U_s$ is open in $B_{s, j}$. Since we also have $U_s \cap B_{s, k} = \emptyset$ for every $k \in I_s \setminus \set{j}$ and $X$ is a tree sum of $\tuple{B_{s, j}: j \in I_s}$ by Observation~\ref{thm:branches}, we have proved the claim.
			
			We also have that $s \in \clo{U_s}$ since $s$ is not isolated in $B_{s, j}$. Finally, let $U := \bigcup\set{U_s: s \in S \cap X_i}$. We have $S \cap X_i \subseteq \clo{U} \setminus U$, and the latter is a closed discrete subset of $X$ by an equivalent condition in Definition~\ref{def:weaker_properties}.
		\end{proof}
	\end{proposition}
	
	\begin{theorem} \labelblock{thm:maximal_connected_tree_sum_equivalence}
		Let $X$ be a tree sum of nondegenerate spaces $\tuple{X_i: i \in I}$. The following conditions are equivalent.
		\begin{enumerate}
			\item The space $X$ is maximal connected.  \loclabel{maximal_connected}
			\item The spaces $X_i$ are maximal connected and $S_X$ is an I-subset of $X$.  \loclabel{I-subset}
			\item The spaces $X_i$ are maximal connected, $S_{X_i}$ is an I-subset of $X_i$ for every $i \in I$, and the gluing is $T_\half$-compatible or, equivalently, I-compatible.  \loclabel{piecewise_I-subset}
			\item The spaces $X_i$ are maximal connected and $X$ is essentially connected.  \loclabel{essentially_connected}
		\end{enumerate}
		
		\begin{proof} \hfill
			\begin{description}
				\item[\locimpl{maximal_connected}{I-subset}.]
					Every $X_i$ is a connected subspace of $X$ by Proposition~\ref{thm:connected_tree_sum}, and hence is maximal connected by Proposition~\itemref{thm:subspace_preservation}{maximal_connected}. Let $S$ denote the set of all non-isolated gluing points of $X$. For every $s \in S$ and $i \in I_s$ we have that $s$ is closed in $X$ since $X$ is $T_\half$. Hence, $s$ is non-isolated in $X_i$ since otherwise it would be clopen in connected nondegenerate space $X_i$. Therefore, $\set{x}$ is nowhere dense in $X_i$, and we may use Proposition~\ref{thm:closed_discrete_gluing_points} to show that $S$ is closed discrete in $X$, and hence $S_X$ is an I-subset of $X$.
				
				\item[\locimpl{I-subset}{maximal_connected}] is Theorem~\itemref{thm:maximal_connected_tree_sum}{maximal_connected}.
				
				\item[\locequiv{I-subset}{piecewise_I-subset}.]
					The equivalence of $T_\half$-compatibility and I-compatibility follows from Observation~\ref{thm:T_half/I-compatibility}. Therefore, the claim follows from Proposition~\ref{thm:gluing_set_discreteness}.
				
				\item[\locequiv{maximal_connected}{essentially_connected}.]
					One implication is trivial, the other follows from Observation~\ref{thm:essential_implies_maximal}.
				\qedhere
			\end{description}
		\end{proof}
	\end{theorem}
	
	\begin{corollary}
		Let $X$ be a tree sum of nondegenerate spaces $\tuple{X_i: i \in I}$. The space $X$ is $T_1$ maximal connected if and only if the spaces $X_i$ are $T_1$ maximal connected and the gluing set $S_X$ is closed discrete.
		
		\begin{proof}
			Follows immediately from Proposition~\ref{thm:tree_sum_separation}, Theorem~\ref{thm:maximal_connected_tree_sum_equivalence}, and the fact that in a nondegenerate connected $T_1$ space there are no isolated points, so every I-subset is closed discrete.
		\end{proof}
	\end{corollary}
\end{subsection}

\end{section}

\begin{section}{Finitely generated spaces}
	
	Maximal connected spaces in the class of finitely generated spaces were first characterized by Thomas in \cite[Theorem 5]{Thomas_68}. He also proposed a way to visualize them. Later, Kennedy and McCartan in \cite{KM_01} characterized finitely generated maximal connected topologies in the lattice $\T(X)$ as joins of two topologies of special form based on the notion of a \emph{final $A$-degenerate cover} where $A \subseteq X$.
	
	In this section we reformulate the characterization in the language of specialization preorder and graphs – finitely generated maximal connected spaces correspond to tree graphs having a fixed bipartition where the correspondence is derived from the graphs of their specialization preorders (Proposition~\ref{thm:finitely_generated_graph_components} and Corollary~\ref{thm:finitely_generated_maximal_connected_graphs}). We also reformulate the characterization in the language of tree sums – they are exactly $T_\half$-compatible tree sums of copies of the Sierpiński space (Corollary~\ref{thm:finitely_generated_maximal_connected_tree_sums}), and we propose another method for their visualization.
	
	\begin{notation}
		Let $X$ be a topological space and $x \in X$. We put 
		\[ \textstyle
			\mon{x} := \bigcap\set{U \subseteq X: U\text{ open neighborhood of }x}.
		\]
	\end{notation}
	
	\begin{observation}
		Let $X$ be a topological space.
		\begin{itemize}
			\item Let $x \in X$. The set $\mon{x}$ is the only candidate for a minimal neighborhood of $x$. Hence, $x$ has a minimal neighborhood if and only if $x^\circ$ is open.
			\item For every, $x, y \in X$ we have $\clo{\set{x}} \subseteq \clo{\set{y}} \iff x \in \clo{\set{y}} \iff \mon{x} \owns y \iff \mon{x} \supseteq \mon{y}$.
			\item $X$ is $T_1$ if and only if $\clo{\set{x}} = \set{x}$ for every $x \in X$ if and only if $\mon{x} = \set{x}$ for every $x \in X$.
			\item $X$ is symmetric if and only if $\clo{\set{x}} = \mon{x}$ for every $x \in X$.
			\item For every $x \in X$ the sets $\clo{\set{x}}$ and $\mon{x}$ are connected.
		\end{itemize}
	\end{observation}

	\begin{definition}
		Recall that a topological space $X$ is called \emph{finitely generated} (or \emph{Alexandrov discrete}) if for every $x \in \clo{A}$ in $X$ there exists a finite set $F \subseteq A$ such that $x \in \clo{F}$, equivalently if every intersection of open sets is open, equivalently if every point $x$ has a minimal neighborhood (which is $\mon{x}$), equivalently if $X$ is inductively generated by the family of all finite subspaces.
	\end{definition}
	
	\begin{observation}
		Every finitely generated $T_\half$ space is submaximal. Hence, a finitely generated space is $T_\half$ if and only if it is submaximal.
		
		\begin{proof}
			Let $X$ be a finitely generated $T_\half$ space. If $D \subseteq X$ dense, then it contains all isolated points, so $\set{x}$ is closed for each $x \in X \setminus D$ because $X$ is $T_\half$, and $X \setminus D$ is closed because $X$ is finitely generated. The last claim follows from Proposition~\ref{thm:implications}.
		\end{proof}
	\end{observation}
	
	\begin{observation} \label{thm:finitely_generated_is_locally_connected}
		Every finitely generated space $X$ is locally connected. In particular, components of connectedness are exactly nonempty clopen connected subsets.
		
		\begin{proof}
			If $x \in U \subseteq X$ for some $U$ open, we have $x \in \mon{x} \subseteq U$ and $\mon{x}$ is open connected.
		\end{proof}
	\end{observation}
	
	\begin{definition} \label{def:specialization_preorder}
		Recall that for every topological space $X$ the \emph{specialization preorder} is defined on its points by the formula 
		\[
			x \leq y \letiff \clo{\set{x}} \subseteq \clo{\set{y}} \iff \mon{x} \supseteq \mon{y}.
		\]
	\end{definition}
	
	The following proposition lists some well-known properties of the specialization preorder.
	We include a proof for the sake of completeness.
	
	\begin{proposition} \labelblock{thm:specialization_preorder}
		Let $X$ be a topological space, $\leq$ the specialization preorder on $X$.
		\begin{enumerate}
			\item Every open set is an upper set. Every closed set is a lower set. \loclabel{open_is_upper}
			\item The converse of \locref{open_is_upper} holds precisely for finitely generated spaces. \loclabel{upper_is_open}
			\item The construction of the specialization preorder provides a $1 : 1$ correspondence between finitely generated spaces and preorders. \loclabel{correspondence}
			\item The specialization preorder is an order if and only if $X$ is $T_0$. \loclabel{order}
			\item Every isolated point is a maximal element. Every closed point is a minimal element. If $X$ is $T_0$ the converse also holds. \loclabel{isolated}
			\item $X$ is $T_\half$ if and only if $\leq$ is an order with at most two levels. \loclabel{T_half}
		\end{enumerate}
		
		\begin{proof} \hfill
			\begin{enumerate}
				\item If $F \subseteq X$ is closed and $x \leq y \in F$, then $x \in \clo{\set{y}} \subseteq F$. Dually for an open set.
				\item An intersection of open sets is an upper set as an intersection of upper sets, and hence it is open if upper sets are open. If $X$ is finitely generated and $U \subseteq X$ is an upper set, then we have $x \in \mon{x} \subseteq U$ for every $x \in U$, and hence $U$ is open. Dually for closed sets.
				\item By \locref{upper_is_open} we know that the construction is injective. We need to show that it is surjective, i.e.\ for every preordered set $\tuple{X, \leq}$ there is a finitely generated topology on $X$ such that $\leq$ is its specialization preorder. It is enough to consider the set of all $\leq$-upper sets as the desired topology. Then $y \in \mon{x}$ if and only if $y$ is in the principal upper set generated by $x$, i.e.\ $y \geq x$.
				\item $\leq$ is an order if and only if $x \in \clo{\set{y}}$ and $y \in \clo{\set{x}}$ implies $x = y$ for every $x, y \in X$.
				\item If $x$ is a closed point, then $y \leq x \iff y \in \clo{\set{x}} = \set{x} \implies y = x$, and hence it is minimal. If $X$ is $T_0$, then $\leq$ is an order, and if $x$ is minimal, then $y \in \clo{\set{x}} \iff y \leq x \implies y = x$, and hence $\set{x}$ is closed. Dually for an isolated point and a maximal element.
				\item If $X$ is $T_\half$, then $\leq$ is an order by \locref{order} and it has at most two levels by \locref{isolated}. If $\leq$ is an order with at most two levels, then $X$ is $T_0$ by \locref{order}, and every point is isolated or closed by \locref{isolated}.
				\qedhere
			\end{enumerate}
		\end{proof}
	\end{proposition}
	
	\begin{definition} \label{def:specialization_graph}
		Let $X$ be a finitely generated $T_\half$ space and $\leq$ its specialization preorder. We define its \emph{specialization graph} $\SG{X}$ as follows. Vertices are the points of $X$, and there is a directed edge $\tuple{x, y}$ in the graph if and only if $x < y$.
	\end{definition}
	
	\begin{proposition} \label{thm:finitely_generated_graph_correspondence}
		The map $X \mapsto \SG{X}$ provides a $1 : 1$ correspondence between finitely generated $T_\half$ spaces and directed graphs with oriented paths of length at most one. On a fixed base set, finer topologies correspond to graphs with less edges.
		
		\begin{proof}
			Clearly, the map $X \mapsto \SG{X}$ factorizes through the construction of specialization preorder, and we have the correspondence between finitely generated $T_\half$ spaces and orders with at most two levels by Proposition \itemref{thm:specialization_preorder}{correspondence}, \itemonlyref{thm:specialization_preorder}{T_half}. Hence, it is enough to establish the correspondence between orders with at most two levels and directed graphs with directed paths of length at most one.
			
			A directed edge $\tuple{x, y}$ is in $\SG{X}$ if and only if $x$ is a closed point, $y$ is an isolated point, and $x \in \clo{\set{y}}$. In that case, $x$ is $\leq$-minimal and $y$ is $\leq$-maximal, and clearly there cannot be oriented paths of length $> 1$ in $\SG{X}$.
			
			On the other hand, we may start with an arbitrary directed graph $G$ with oriented paths of length at most one and interpret it as a strict part of an order with at most two levels.
		\end{proof}
	\end{proposition}
	
	\begin{lemma} \label{thm:simple_expansion_of_finitely_generated}
		Every simple expansion of a finitely generated space is finitely generated. Hence, every simple expansion of a finitely generated $T_\half$ space is finitely generated and $T_\half$.
		
		\begin{proof}
			Let $\tuple{X, \tau}$ be a finitely generated topological space, let $A \subseteq X$, and $\tau^* := \tau \vee \set{A}$. Every family $\U$ of $\tau^*$-open sets is of form $\set{(U_i \cup A) \cap V_i: i \in I}$ where all the sets $U_i$, $V_i$ are $\tau$-open. Hence, $\bigcap\U = \bigcap_{i \in I} (U_i \cup A) \cap V_i = ((\bigcap_{i \in I} U_i) \cup A) \cap (\bigcap_{i \in I} V_i)$ is $\tau$-open.
		\end{proof}
	\end{lemma}
	
	\begin{proposition} \labelblock{thm:finitely_generated_graph_components}
		Let $X$ be a finitely generated $T_\half$ space.
		\begin{enumerate}
			\item Connected components of $X$ are exactly undirected connected components of $\SG{X}$. Hence, $X$ is connected if and only if $\SG{X}$ is connected (as undirected graph).\loclabel{connected}
			\item $X$ is maximal connected if and only if $\SG{X}$ is a tree (as undirected graph). \loclabel{maximal_connected}
		\end{enumerate}
		
		\begin{proof} \hfill
			\begin{enumerate}
				\item Let $x \in X$. The connected component of $X$ containing $x$ is the lower upper set generated by $\set{x}$ because that is the smallest clopen set containing $x$ (we use Proposition \ref{thm:finitely_generated_is_locally_connected}). That is exactly the component of undirected connectedness of $\SG{X}$.
				\item We use the correspondence from Proposition~\ref{thm:finitely_generated_graph_correspondence}. Connected topologies correspond to connected graphs, and trees are exactly the minimal connected graphs. We also need Lemma~\ref{thm:simple_expansion_of_finitely_generated} and the fact that maximal connectedness can be tested on simple expansions (Observation~\ref{thm:testing_expansions}).
				\qedhere
			\end{enumerate}
		\end{proof}
	\end{proposition}
	
	\begin{corollary} \label{thm:finitely_generated_maximal_connected_graphs}
		Finitely generated maximal connected spaces correspond to tree graphs with a fixed bipartition.
		
		\begin{proof}
			We can use Proposition~\itemref{thm:finitely_generated_graph_components}{maximal_connected} because every maximal connected space is $T_\half$ by Proposition~\ref{thm:implications}, and because we can equivalently describe directed graphs with directed paths of length at most one as undirected graphs with a fixed bipartition.
		\end{proof}
	\end{corollary}
	
	\begin{proposition} \label{thm:finitely_generated_tree_sum}
		Let $X$ be a tree sum of spaces $\tuple{X_i: i \in I}$. The space $X$ is finitely generated if and only if all spaces $X_i$ are finitely generated.
		
		\begin{proof}
			One implication is clear since every subspace of a finitely generated space is itself finitely generated. On the other hand, if every $X_i$ is finitely generated, it is inductively generated by its finite subspaces, and since $X$ is inductively generated by the spaces $X_i$, it is, by transitivity, inductively generated by some of its finite subspaces and hence by all of its finite subspaces.
		\end{proof}
	\end{proposition}
	
	\begin{corollary} \label{thm:finitely_generated_maximal_connected_tree_sums}
		Besides the one-point space, finitely generated maximal connected spaces are exactly $T_\half$-compatible tree sums of copies of the Sierpiński space.
		
		\begin{proof}
			Clearly, the Sierpiński space is maximal connected. Hence, every $T_\half$-compatible tree sum of copies of the Sierpiński space is maximal connected by Theorem~\ref{thm:maximal_connected_tree_sum_equivalence}, and it is finitely generated by Proposition~\ref{thm:finitely_generated_tree_sum}.
			
			On the other hand, let $X$ be a finitely generated maximal connected space. Clearly, every edge of $\SG{X}$ corresponds to a subspace of $X$ homeomorphic to the Sierpiński space. By Proposition~\ref{thm:finitely_generated_graph_components} the corresponding gluing graph is a tree, and the subspaces cover whole $X$ unless $X$ is a one-point space. Also, the Sierpiński subspaces inductively generate the topology. If $x \in \clo{A} \setminus A$ for some $A \subseteq X$, then there is $y \in A$ such that $x \in \clo{\set{y}}$, so $x < y$ and $\set{x, y}$ is the witnessing Sierpiński subspace. Altogether, $X$ a tree sum of the Sierpiński subspaces by Proposition~\ref{thm:inner_tree_sum}, and the gluing is $T_\half$-compatible again by Theorem~\ref{thm:maximal_connected_tree_sum_equivalence}.
		\end{proof}
	\end{corollary}
	
	\begin{notation}
		We propose to visualize finitely generated maximal connected spaces as follows. We just draw the specialization graph, and instead of orienting the edges we distinguish two kinds of vertices: the vertices corresponding to isolated points shall be drawn as open dots while the vertices corresponding to closed points as solid dots. Whenever it is suitable, we draw the isolated vertices above the closed vertices in order to stress the specialization order.
		
		In the original paper \cite{Thomas_68} Thomas used a similar but different visualization. Instead of drawing an edge from a closed point to every isolated point in its smallest neighborhood, Thomas represents the smallest neighborhood by a line segment containing all the points.
		
		We think our visualization is less restrictive, and it deals with the duality on finitely generated spaces well – it is enough to switch the colors of isolated and closed vertices.
	\end{notation}
	
	\begin{example}
		We give some examples of maximal connected finitely generated spaces. Visualizations of these spaces are given in Figure~\ref{fig:spaces}.
		\begin{itemize}
			\item Clearly, the empty space, the one-point space, and the Sierpiński space are maximal connected.
			\item If $X$ is a set, $x \in X$ and $A \subseteq X$ is open if and only if $x \in A$, we obtain a space with so-called \emph{included point topology}, also called \emph{principal ultrafilter space}.
			\item If $X$ is a set, $x \in X$ and $A \subseteq X$ is open if and only if $x \notin A$, we obtain a space with so-called \emph{excluded point topology}, also called \emph{principal ultraideal space}.
			\item Let us consider the set of all integers $\ZZ$ with the topology generated by open sets $\set{\set{2k - 1,\, 2k,\, 2k + 1}: k \in \ZZ}$. We obtain a finitely generated maximal connected space called \emph{Khalimsky line} or \emph{digital line}.
		\end{itemize}
	\end{example}
	
	\input{figures/spaces.fig}
	
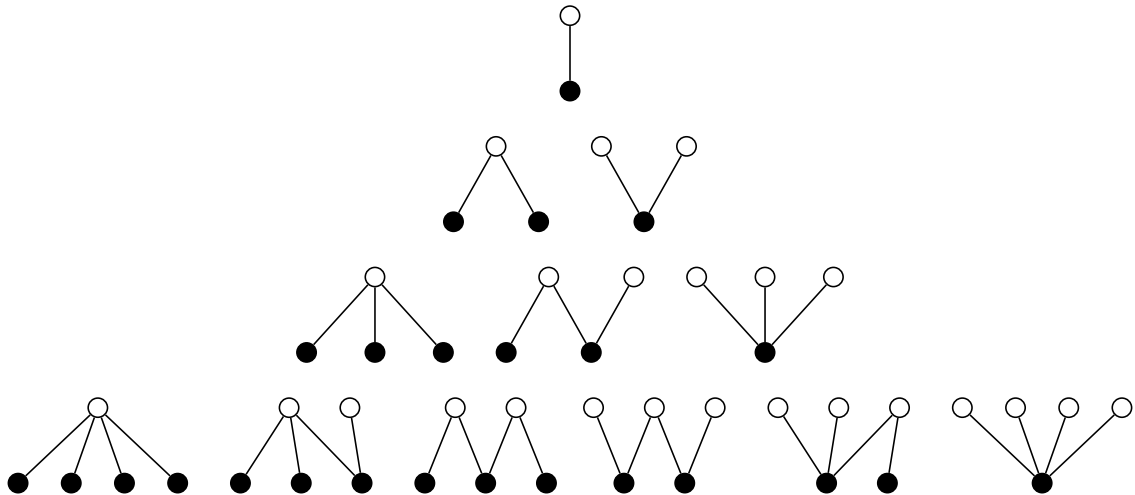
	Let us conclude with Figure~\ref{fig:small_spaces}, which shows all nondegenerate finitely generated maximal connected spaces with at most five points, using our visualization. This may be compared with the corresponding picture in \cite{Thomas_68}. In our picture, the duality of finitely generated spaces (by considering the dual specialization preorder, in our case just by swapping isolated and closed points) is apparent.
	
	\input{figures/small_spaces.fig}
\end{section}

	\begin{acknowledgements}
		The author was supported by the grant SVV-2016-260336 of Charles University and by the University Center for Mathematical Modeling, Applied Analysis and Computational Mathematics (MathMAC).
		The author would also like to thank the referee for his or her comments.
	\end{acknowledgements}
	
	\linespread{1}\selectfont
	
	\bibliographystyle{siam}
	\bibliography{references}

\begin{thebibliography}{10}

\bibitem{Anderson_65}
{\sc D.~R. Anderson}, {\em On connected irresolvable {H}ausdorff spaces}, Proc.
  Amer. Math. Soc., 16 (1965), pp.~463--466.

\bibitem{AC_95}
{\sc A.~V. Arhangel'ski{\u\i} and P.~J. Collins}, {\em On submaximal spaces},
  Topology Appl., 64 (1995), pp.~219--241.

\bibitem{Cameron_71}
{\sc D.~E. Cameron}, {\em Maximal and minimal topologies}, Trans. Amer. Math.
  Soc., 160 (1971), pp.~229--248.

\bibitem{Engelking_89}
{\sc R.~Engelking}, {\em General topology}, vol.~6 of Sigma Series in Pure
  Mathematics, Heldermann Verlag, Berlin, second~ed., 1989.
\newblock Translated from the Polish by the author.

\bibitem{GRS_73}
{\sc J.~A. Guthrie, D.~F. Reynolds, and H.~E. Stone}, {\em Connected expansions
  of topologies}, Bull. Austral. Math. Soc., 9 (1973), pp.~259--265.

\bibitem{GS_73}
{\sc J.~A. Guthrie and H.~E. Stone}, {\em Spaces whose connected expansions
  preserve connected subsets}, Fund. Math., 80 (1973), pp.~91--100.

\bibitem{GSW_78}
{\sc J.~A. Guthrie, H.~E. Stone, and M.~L. Wage}, {\em Maximal connected
  expansions of the reals}, Proc. Amer. Math. Soc., 69 (1978), pp.~159--165.

\bibitem{Hildebrand_67}
{\sc S.~K. Hildebrand}, {\em A connected topology for the unit interval}, Fund.
  Math., 61 (1967), pp.~133--140.

\bibitem{KM_01}
{\sc G.~J. Kennedy and S.~D. McCartan}, {\em Maximal connected principal
  topologies}, Math. Proc. R. Ir. Acad., 101A (2001), pp.~163--166.

\bibitem{Levine_70}
{\sc N.~Levine}, {\em Generalized closed sets in topology}, Rend. Circ. Mat.
  Palermo (2), 19 (1970), pp.~89--96.

\bibitem{McSherry_74}
{\sc D.~M.~G. McSherry}, {\em On separation axioms weaker than {$T_{1}$}},
  Proc. Roy. Irish Acad. Sect. A, 74 (1974), pp.~115--118.

\bibitem{Simon_78}
{\sc P.~Simon}, {\em An example of maximal connected {H}ausdorff space}, Fund.
  Math., 100 (1978), pp.~157--163.

\bibitem{Thomas_68}
{\sc J.~P. Thomas}, {\em Maximal connected topologies}, J. Austral. Math. Soc.,
  8 (1968), pp.~700--705.

\bibitem{van_Douwen_93}
{\sc E.~K. van Douwen}, {\em Applications of maximal topologies}, Topology
  Appl., 51 (1993), pp.~125--139.

\end{thebibliography}
	
\end{document}